\documentclass[12pt]{article}

\usepackage[T2A]{fontenc}
\usepackage[cp1251]{inputenc}
\usepackage[russian]{babel}
\usepackage{amsfonts}
\newtheorem{thm}{Теорема}[section]
\newtheorem{lemma}{Лемма}[section]
\newtheorem{conseq}{Следствие}[section]
\newtheorem{definition}{Определение}
\newtheorem{condition}{Условие}[section]
\newtheorem{agreement}{Соглашение}

\begin{document}

\title{Равномерная оценка коэффициента сильного перемешивания и максимума остаточного эмпирического процесса в ARCH($p$) модели.}
\author{А. А. Сорокин}
\date{}
\maketitle

\section{Введение.}

На протяжении нескольких последних десятилетий гетероскедастическим моделям было посвящено значительное число теоретических исследований. Они также находят широкое применение на практике. Базовой для всего семейства гетероскедастических моделей является ARCH($p$) модель, предложенная Робертом Энглем (Robert F. Engle) в 1982 году в работе \cite{Engle}.

\begin{definition}
ARCH($p$) модель определяется как решение следующего уравнения:
\begin{equation}\label{equ_arch_p_prel}
y_t = \sigma_t(\mathbf{a}) \varepsilon_t, \quad \sigma^2_t(\mathbf{a}) = a_0 + a_1 y_{t-1}^2 + \ldots + a_p y_{t-p}^2, \quad t \in \mathbb{Z},
\end{equation}
где $\{y_t\}$-наблюдения (за логарифмическими приращениями цены некоторого актива), $\mathbf{a} := (a_0, a_1, \ldots, a_p)^*$ - вектор неизвестных параметров ($*$ - знак транспонирования),  $a_0 > 0, a_1 \geq 0, \ldots, a_p \geq 0$, $\sigma_t(\mathbf{a})$ - функция волатильности, $\{\varepsilon_t\}$ - н.о.р.с.в. с неизвестной функцией распределения $G(x).$
\end{definition}

После появления ARCH($p$) модели было предложено множество ее изменений и обобщений. Мы упомянем только наиболее часто цитируемую и употребляемую на практике GARCH($p, q$) модель (см. \cite{Bollerslev}).

Одной из самых интересных и важных статистических задач для ARCH($p$) модели является оценка вектора параметров $\mathbf{a}$. Наиболее известными типами оценок являются оценки максимального правдоподобия (см., например, обзор \cite{Shephard}) и оценки квази-максимального правдоподобия (см. \cite{LH}, \cite{Rahbek}, \cite{Straumann}).

В работах \cite{Boldin}, \cite{Boldin3}, \cite{S}, \cite{S2}, \cite{Vyazilov}, предлагался метод исследования оценок и тестов для (G)ARCH модели, основанный на понятии остаточного эмпирического процесса. Он дает возможность единым способом доказать асимптотическую нормальность различных типов оценок, используя линейное разложение остаточного эмпирического процесса. Такой метод также хорошо известен для линейных моделей (см., например, \cite{Koul}, \cite{Koul2}).

Прежде чем перейти к построению остаточного эмпирического процесса и формулировке основных результатов работы, заметим, что модель (\ref{equ_arch_p_prel}) с неизвестной $G(x)$ не является \emph{идентифицируемой}. Действительно, замена вида $\varepsilon_t' = \varepsilon_t / \gamma, \; a_0' = a_0 \gamma^2, \ldots, a_p' = a_p \gamma^2$ с $\gamma > 0$ приводит к тому же процессу $\{y_t' \equiv y_t\}.$ Мы положим $a_0 = 1$, тем самым подразумевая, что наблюдаемый процесс $\{y_t\}$ является реализацией модели (\ref{equ_arch_p_prel}) с $a_0 = 1.$ Для этого достаточно сделать описанную замену с $\gamma := (a_0)^{-1/2}$. Таким образом, в дальнейшем мы рассматриваем модель
\begin{equation}\label{equ_arch_p}
y_t = \sigma_t(\mathbf{a}) \varepsilon_t, \quad \sigma^2_t(\mathbf{a}) = 1 + a_1 y_{t-1}^2 + \ldots + a_p y_{t-p}^2, \quad t \in \mathbb{Z},
\end{equation}
Мы будем опускать $a_0$ в записи $\mathbf{a}$ и обозначать $\mathbf{a} := (a_1, \ldots, a_p)^*$. Предположим, что $G(x)$ абсолютно непрерывна относительно меры Лебега с плотностью $g(x)$. Пусть выполнено следующее
\begin{condition}\label{cond_station}
$$
E \varepsilon_1 = 0, \; E \varepsilon_1^2 < \infty, \; E \varepsilon_1^2 (a_1 + \ldots + a_p) < 1, \; g(x) > 0.
$$
\end{condition}
Тогда (см. \cite[стр. 106-107]{Doukhan}), уравнение (\ref{equ_arch_p}) имеет единственное строго стационарное решение, обозначаемое в дальнейшем $\{y^a_t\}$ или просто $\{y_t\}$. Будем также полагать $\mathbf{Y}^a_t = (y^a_t, \ldots, y^a_{t - p + 1})^*$.
\begin{agreement}
Обычно для сокращения записи мы будем опускать верхний индекс $\mathbf{a}$ у величин, зависящих от $\{y^a_t\}$. Например, будем полагать $y_t = y_t^a, \; \mathbf{Y}_t = \mathbf{Y}^a_t, \ldots$.
\end{agreement}
С целью сокращения изложения будем предполагать, что $E \varepsilon_1^2 \geq \beta_a$ для некоторого заранее известного $\beta_a > 0$ (основные результаты могут быть доказаны и без такого предположения). Пусть $\|\cdot\|_1$ и $\|\cdot\|$ обозначают соответственно $L_1$ и $L_2$ нормы. С учетом условия \ref{cond_station}, очевидно, имеем $\|\mathbf{a}\|_1 < \beta_a^{-1}.$ Таким образом, область допустимых значений для $\mathbf{a}$ имеет вид
$$
\Theta := \{\mathbf{a} \in \mathbb{R}^p: a_1 \geq 0, \ldots a_p \geq 0, \; \|\mathbf{a}\|_1 \leq \beta_a^{-1}\}.
$$

Перейдем к построению остаточного эмпирического процесса для модели (\ref{equ_arch_p}). Положим для $\theta \in \Theta, \; \mathbf{U} \in \mathbb{R}^{p}, \; t = 1, \ldots, n$
$$
s(\theta, \mathbf{U}) :=  1 + \theta_1 U_1^2 + \ldots + \theta_p U_p^2, \quad \sigma^2_t(\theta) := s(\theta, \mathbf{Y}_{t-1}).
$$
Определим остатки для уравнения (\ref{equ_arch_p}) соотношением $\varepsilon_t(\theta) := y_{t} \sigma^{-1}_t(\theta), \; \theta \in \Theta$. Положим
$$
G_n(x, \theta) := n^{-1} \sum \limits_{t=1}^n I\{\varepsilon_t(\theta) \leq x\},
$$
где $I\{\cdot\}$ - индикатор события. Функция $G_n(x, \theta)$ называется остаточной эмпирической ф.р. Пусть $\varphi(\mathbf{U}, \theta) := (\varphi_1(\mathbf{U}, \theta), \ldots, \varphi_p(\mathbf{U}, \theta))^*-$ некоторая функция, $\varphi : \mathbb{R}^p \times \Theta \to \mathbb{R}^p$.

\begin{definition}
Остаточный эмпирический процесс $\mathbf{W}_n(x, \theta)$ определяется соотношениями
$$
\mathbf{W}_n(x, \theta) := (W_{n, 1}(x, \theta), \ldots, W_{n, p}(x, \theta))^*,
$$
$$
W_{n, j}(x, \theta) := n^{-1/2} \sum \varphi_j(\mathbf{Y}_{t-1}, \theta) [ I\{\varepsilon_t(\theta) \leq x\} - G_n(x, \theta)], j = 1, \ldots, p,
$$
где $n \in \mathbb{N}$, $x \in \mathbb{R}$, $\theta \in \Theta$.
\end{definition}
\begin{agreement}
Здесь и в дальнейшем множество интегрирования совпадает с $\mathbb{R}$, суммирование проводится по $t$ от $1$ до $n$, и предельный переход осуществляется при $n \to \infty$, если явно не указано обратное.
\end{agreement}

В работах \cite{Boldin}, \cite{Vyazilov}, \cite{S}, \cite{S2} рассматривались два типа оценок - оценки минимального расстояния и Generalized M оценки. Напомним их определение.

\begin{definition}
Пусть фиксирована некоторая ф.р. $F(x)$. Тогда решение экстремальной задачи
$$
K_n(\theta) := \int \|\mathbf{W}_n(x, \theta)\| d F(x) \to \min \limits_{\theta \in \Theta}
$$
называется оценкой минимального расстояния (в дальнейшем MD).
\end{definition}

\begin{definition}
Пусть фиксирована функция $\psi(x)$. Тогда решение (векторного) уравнения
$$
\mathbf{l}_n(\theta) := \sum \varphi(\mathbf{Y}_{t-1}, \theta) \psi(\varepsilon_t^n(\theta)) - \sum \varphi(\mathbf{Y}_{t-1}, \theta) \overline{\psi}_n(\theta) = \mathbf{0},
$$
где $\overline{\psi}_n(\theta) := n^{-1} \sum \psi(\varepsilon^n_t(\theta))$, называется Generalized M оценкой (в дальнейшем GM).
\end{definition}

Несложно заметить, что целевые функционалы (левые части) для обоих типов оценок могут быть выражены как функции от о.э.п. Для MD оценки это очевидно, а в случае GM имеет место соотношение $\mathbf{l}_n(\theta) = - \int \mathbf{W}_n(x, \theta) d \psi(x)$. Поэтому для доказательства асимптотической нормальности оценок может быть применена техника, основанная на равномерной асимптотической линейности о.э.п. Введем дополнительные определения. Положим для $\mathbf{U} \in \mathbb{R}^p$, $\theta \in \Theta, \; k = 1, \ldots, p$
$$
e_k(\mathbf{U}, \theta) := U^2_k s(\theta, \mathbf{U})^{-1}, \quad \mathbf{e}^{0}(\mathbf{U}, \theta) := (e^{0}_1(\mathbf{U}, \theta), \ldots, e^{0}_p(\mathbf{U}, \theta))^*,
$$
$$
e^{0 a}_k(\mathbf{U}, \theta) := e_k(\mathbf{U}, \theta) - E e_k(\mathbf{Y}^a_0, \theta), \quad \mathbf{e}^{0 a}(\mathbf{U}, \theta) := (e^{0 a}_1(\mathbf{U}, \theta), \ldots, e^{0 a}_p(\mathbf{U}, \theta))^*,
$$
$$
\varphi^{0 a}_k(\mathbf{U}, \theta) :=  \varphi_k(\mathbf{U}, \theta) - E \varphi_k(\mathbf{Y}^a_0, \theta), \quad \varphi^{0 a}(\mathbf{U}, \theta) := (\varphi^{0 a}_1(\mathbf{U}, \theta), \ldots, \varphi^{0 a}_p(\mathbf{U}, \theta))^*,
$$
$$
S_{\varphi, e}(\mathbf{a}) := E \varphi^{0 a}(\mathbf{Y}^a_0, \mathbf{a}) [\mathbf{e}^{0 a}( \mathbf{Y}^a_0, \mathbf{a})]^*, \quad S_{\varphi, \varphi}(\mathbf{a}) := E \varphi^{0 a}(\mathbf{Y}^a_0, \mathbf{a}) [\varphi^{0 a}(\mathbf{Y}^a_0, \mathbf{a})]^*,
$$
$$
\mathbf{\widetilde{W}}_n(x, \theta) := \mathbf{W}_{n}(x, \mathbf{a}) + (1/2) n^{1/2} x g(x)  S_{\varphi, e}(\mathbf{a}) (\theta - \mathbf{a}).
$$
\begin{definition}[См. \cite{Koul}, Глава 5.]
Будем говорить, что для о.э.п. выполнено свойство равномерной асимптотической линейности (Asympto\-tic Uni\-form Line\-arity, AUL), если для любого $B < \infty$
$$
\sup \limits_{\|s\| \leq B} \|\mathbf{\widetilde{W}}_n(x, \mathbf{a} + n^{-1/2} \mathbf{s}) - \mathbf{W}_n(x,\mathbf{a} + n^{-1/2} \mathbf{s})\| = o_p(1).
$$
\end{definition}
Как показывается в работах \cite{Koul}, \cite{Koul2}, из AUL легко следует асимптотическая нормальность для широкого класса оценок в линейных по $\varepsilon$ моделях, основанных на о.э.п. (например для MD), в предположении их $n^{1/2}$-состоятельности. Для линейных моделей $n^{1/2}$-состоятельность автоматически имеет место, так как о.э.п. монотонно зависит от $\theta$ (см. детали в \cite{Koul}). В случае ARCH($p$) модели это не так. Для ее установления в работе \cite{S} был предложен следующий метод. Было доказано, что равномерно вне $n^{-1/2}$-окрестности $\mathbf{a}$ о.э.п. принимает "большие" \, значения. С этой целью использовалось специальное максимальное неравенство для о.э.п. Введем еще одно определение. Положим для $\theta \in \Theta, \; \mathbf{a} \in \Theta$
$$
\mathbf{b}^a(x, \theta) := E \mathbf{\varphi}(\mathbf{Y}^a_0, \theta)
I\{\varepsilon^a_1(\theta) \leq x \} - E \mathbf{\varphi}(\mathbf{Y}^a_0, \theta) P(\varepsilon^a_1( \theta) \leq x).
$$
\begin{definition}
Будем говорить, что о.э.п. ограничен по вероятности с точностью до $\|\theta - \mathbf{a}\|$, если для любых $\rho > 0$, $k = 1, \ldots, p$
\begin{equation}\label{equ_wn_rough_estim}
\sup \limits_{\theta \in \Theta} (W_{n, k}(x, \theta) - n^{1/2} b^a_k(x, \theta) \pm n^{1/2} \rho \|\theta - \mathbf{a}\| )^{\mp} = O_p(1).
\end{equation}
\end{definition}
В работе \cite{S} аналог (\ref{equ_wn_rough_estim}) устанавливался для о.э.п. в ARCH($1$) модели. Можно проверить, что следующие условия являются достаточными для $n^{1/2}$-состоятельности соответственно MD и GM оценки, в предположении справедливости (\ref{equ_wn_rough_estim}).
\begin{condition}\label{cond_phi_info_md}
Для некоторого $\delta_0 > 0$ и любых $\mathbf{a} \in \Theta^{\delta}, \; \theta \in \Theta$ $$
\int \| \mathbf{b}^a(x, \theta) \|^2 d F(x) \geq \delta_0 ||\theta - \mathbf{a}||^2.
$$
\end{condition}
\begin{condition}\label{cond_phi_info_gm}
Для некоторого $\delta_0 > 0$ и любых $\mathbf{a} \in \Theta^{\delta}, \; \theta \in \Theta$ $$
\left \| \int \mathbf{b}^a(x, \theta) d \psi(x) \right \|  \geq \delta_0 ||\theta - \mathbf{a}||.
$$
\end{condition}

Как при доказательстве AUL, так и при доказательстве ограниченности по вероятности с точностью до $\|\theta - \mathbf{a}\|$ о.э.п., важную роль играет оценка скорости убывания зависимости процесса $\{y_t^a\}$. В первом случае достаточно просто эргодичности (см., например, см. \cite{Koul}, \cite{Koul2}, \cite{Boldin}). Для доказательства аналога (\ref{equ_wn_rough_estim}) необходимо более сильное ограничение. В работе \cite{S} (см. \cite[Corollary 5.1]{S}) был использован тот факт, что коэффициент сильного перемешивания процесса $\{y_t^a\}$ убывает экспоненциально быстро. Его доказательство для случая ARCH($p$) модели можно найти, например, в \cite{Doukhan}.

Геометрическая оценка коэффициента сильного перемешивания также была обобщена несколькими авторами на GARCH($p, q$) модель. Соответствующий факт содержится, например, в работах \cite{Boussama} (на французском) и \cite{SM2}, \cite{Straumann}. В двух последних работах доказательство основывается на оценке коэффициента сильного перемешивания процесса, являющегося решением так называемого Linear Polynomial Stochastic Recurrence Equation (SRE). В наших обозначениях последнее имеет вид
$$
\mathbf{Y}_t = \mathbf{P}_t(\varepsilon_t) \mathbf{Y}_{t - 1} + \mathbf{Q}_t(\varepsilon_t),
$$
где $\mathbf{P}_t(x), \; \mathbf{Q}_t(x)$ - матрицы, каждый элемент которых является полиномом от $x$. Общий результат об оценке коэффициента сильного перемешивания для SRE содержится в работе \cite{Mokkadem} (на французском) и, в немного измененном виде, он цитируется в \cite{SM2} (Теорема 4.5).

В настоящей работе устанавливаются два основных результата. Первой из них является обобщением предыдущего выше утверждения про коэффициент сильного перемешивания. Второй устанавливает оценку для максимума остаточного эмпирического процесса общего вида. Прежде всего определим схему наблюдений, для которой будут формулироваться наши результаты.
\begin{definition}
Пусть фиксированы некоторые $\mathbf{b} \in \Theta$ и $\delta > 0$, для которых $E \varepsilon_1^2 (||\mathbf{b}||_1 + \delta) < 1$. Предположим, что для $n \geq 1$ выборка состоит из величин $y^n_{1-p}, \ldots, y^n_0, y^n_1, \ldots, y^n_n,$ являющихся выборкой из стационарного решения (\ref{equ_arch_p}) с параметром
\begin{equation}\label{equ_observ_scheme}
\mathbf{a}_n \in \Theta^{\delta} := \{\mathbf{a} \in \mathbb{R}^p: \|\mathbf{a} - \mathbf{b}\|_1 \leq \delta \} \cap \Theta,
\end{equation}
зависящим от $n$ и принадлежащим $\Theta^{\delta}$. В дальнейшем мы будем говорить о (\ref{equ_observ_scheme}) как о схеме наблюдений с зависящим от $n$ параметром.
\end{definition}
Чтобы подчеркнуть зависимость параметра от $n$, мы будем обозначать $\{y^{n}_t\}$ решение уравнения (\ref{equ_arch_p}) с $\mathbf{a} = \mathbf{a}_n$. Тем самым $y^{a_n}_t = y^n_t.$ Будем также для $t = 1, \ldots, n$ полагать
$$
\mathbf{Y}^n_{t} = \mathbf{Y}^{a_n}_{t}, \quad \sigma^2_{t n}(\theta) := s(\theta, \mathbf{Y}^n_{t - 1}),
$$
$$
\varepsilon_t^n(\theta) := y_t^n / \sigma_{t n}(\theta), \quad G_n(x, \theta) := n^{-1} \sum I\{\varepsilon_t^n(\theta) \leq x\}.
$$
и определять компоненты о.э.п. как
$$
W_{n, j}(x, \theta) := n^{-1/2} \sum \varphi_j(\mathbf{Y}^n_{t-1}, \theta) [ I\{\varepsilon^n_t(\theta) \leq x\} - G_n(x, \theta)], j = 1, \ldots, p.
$$

Для установления AUL в такой схеме наблюдений мы будем использовать \emph{равномерную} по $\mathbf{a} \in \Theta^{\delta}$ оценку коэффициента сильного перемешивания $\{y^a_t\}$. К сожалению, из результатов \cite{Mokkadem}, \cite{SM2} этот факт прямо не следует. Поэтому мы приводим свое доказательство нужного результата, полученное независимо. Это утверждение (следствие \ref{thm_expan_recur_mixing_coef_arch}) составляет первый из наших результатов.

Второй результат (теорема \ref{thm_general_proc_op}) есть оценка максимума остаточного эмпирического процесса общего вида, построенного по процессу (нелинейной) авторегрессии, обладающего достаточно быстрым сильным перемешиванием. В частности, будет проверено, что оценка имеет место для ARCH($p$) модели (следствие \ref{thm_arch_p_proc_op}). Аналогичные результаты в литературе по гетероскедастическим моделям нам неизвестны.

В качестве дополнительной иллюстрации полученных результатов нами рассматривается следующая задача. Мы устанавливаем робастность GM и MD оценок против грубых засорений (тем самым обобщая теорему 2.2 из \cite{S}). При этом используется максимальное неравенство, устанавливаемое следствием \ref{thm_expan_recur_mixing_coef_arch}.

Схема дальнейшего изложения выглядит следующим образом. В разделе \ref{seq_results} доказываются оценка максимума для о.э.п. в ARCH($p$) модели (параграф \ref{seq_maq_inequality_corrol}) и робастность MD и GM оценок (параграф \ref{seq_robustness}). В разделе \ref{seq_technical_results} доказываются результаты об оценке коэффициента сильного перемешивания для ARCH($p$) модели (параграф \ref{seq_regression_uniform_mixing}) и максимальное неравенство для о.э.п. общего вида (параграф \ref{seq_max_inequality}).

\section{Робастность MD и GM оценок.} \label{seq_results}

Прежде чем перейти к формулировке результатов, введем необходимые технические условия:
\begin{condition}\label{cond_distr}
$$
g(x) > 0, \lim \limits_{|x| \rightarrow \infty} x g(x) = 0, \sup \limits_{x} (1 + x^2) |g'(x)| < \infty, \; \int | g'(x) x| d x < \infty.
$$
\end{condition}
\begin{condition}\label{cond_deriv}
Существуют функции $\Phi(\cdot)$, $M(\cdot)$, а также константа $\beta_0 > 0$, для которых

i) Для $\mathbf{U} \in \mathbb{R}^p, \; \theta \in \Theta$ $|\varphi_j(\mathbf{U}, \theta)| \leq \Phi(\mathbf{U}), \; \theta \in \Theta,$ и
$$
\sup \limits_{a \in \Theta^{\delta}} E \Phi^{4 + \beta_0}(\mathbf{Y}^a_0) < \infty.
$$

ii) Для $\theta \in \Theta^{\delta}$, $j = 1, \ldots, p$ и п.в. $\mathbf{U} \in \mathbb{R}^p$ существует непрерывная по $\theta$ $d_{j, k}( \mathbf{U}, \theta) := \frac{\partial{\varphi_j( \mathbf{U}, \theta)}} {\partial {\theta_k}},$ причем
$$
\|d_{j, k}( \mathbf{U}, \theta)\| \leq M(\mathbf{U}), \quad \sup \limits_{a \in \Theta^{\delta}} E M^{2 + \beta_0}(\mathbf{Y}^a_0) < \infty.
$$
\end{condition}
Нам также дополнительно понадобится следующее ограничение:
\begin{condition}\label{cond_epsilon_mom_8_plus_e}
Для некоторого $\beta_0 > 0$ $E |\varepsilon_1|^{8 + \beta_0} < \infty$.
\end{condition}
Условия \ref{cond_distr} - \ref{cond_epsilon_mom_8_plus_e} являются техническими и, вероятно, могут быть ослаблены.

\subsection{Оценка максимума о.э.п. для ARCH($p$) модели.} \label{seq_maq_inequality_corrol}

В данном параграфе мы установим максимальное неравенство для $\mathbf{W}_n(x, \theta)$. С этой целью будет применена теорема \ref{thm_general_proc_op} настоящей работы.
\begin{agreement}
Будем обозначать $O^{a}_p(1)$ (соответственно $o^{a}_p(1)$) семейство с.в. $\{\eta^{a}_n\}$ , для которого $\lim \limits_{C \to \infty} \sup \limits_{a \in \Theta^{\delta}, n} P(|\eta_n^{a}| \geq C) = 0$ (соответственно $\lim \limits_{n \to \infty} \sup \limits_{a \in \Theta^{\delta}} P(|\eta^{a}_n| \geq \delta) = 0$).
\end{agreement}

\begin{conseq}\label{thm_arch_p_proc_op}
Пусть выполнены условия \ref{cond_station}, \ref{cond_distr} - \ref{cond_epsilon_mom_8_plus_e}. Тогда для $k = 1, \ldots, p$ и любого $\rho > 0$
\begin{equation}\label{equ_wn_exp_wn_diff_estim}
\sup \limits_{x \in \mathbb{R}, \theta \in \Theta} (W_{n, k}(x, \theta) - n^{1/2} [b^{a_n}_k(x, \theta) \pm \rho \| \theta - \mathbf{a}_n \|])^{\pm} = O^a_p(1),
\end{equation}
\end{conseq}
\textbf{Доказательство.} Положим $l^a(x, \theta, \mathbf{U}) := x s(\theta, \mathbf{U})^{1/2} s(\mathbf{a}, \mathbf{U})^{-1/2}$, таким образом $I\{\varepsilon^a_t(\theta) \leq x\} \equiv I\{\varepsilon_t \leq l^a(x, \theta, \mathbf{Y}^a_{t-1})\}$. Несложно проверить, что
$$
2^{-1} x g(x) e_k(\mathbf{U}, \mathbf{a}) = \left. \frac{\partial{G(l^a(x, \theta, \mathbf{U} ))}}{\partial{\theta_k}} \right |_{\theta = a}.
$$
Обозначим для краткости при $k = 1, \ldots, p$ $\varphi^n_{k, t}(\theta) := \varphi_k(\theta, \mathbf{Y}^a_{t-1}),$ $v^a_t(x, \theta) := I\{\varepsilon^a_t(\theta) \leq x\},$ в таком случае для $k = 1, \ldots, p$
$$
b^{a}_k(x, \theta) = [E \varphi^n_{k, 1}(\theta) v^n_1(x, \theta) - E \varphi^n_{k, 1}(\theta) E v^n_1(x, \theta)].
$$
Для доказательства (\ref{equ_wn_exp_wn_diff_estim}) проверим, что для $k = 1, \ldots, p$ выполнены условия теоремы \ref{thm_general_proc_op} с
$$
q = p, \; \beta = \beta_0 / 2, \; \mathbf{z} = \theta - \mathbf{a}, \; \tau = \mathbf{a}, \; \mathbf{z}_0 = \mathbf{0},
$$
$$
\lambda^{\tau}(\mathbf{z}, \mathbf{U}) = \varphi_k(\theta, \mathbf{U}), \; \Delta^{\tau}(x, \mathbf{z}, \mathbf{U}) = l^a(x, \theta, \mathbf{U}).
$$
Условие i) выполняется автоматически в силу определения $l^a(x, \theta, \mathbf{U}).$ Следствие \ref{thm_expan_recur_mixing_coef_arch} влечет за собой справедливость условия iv). Проверим условие ii). Положим $g_R := \sup \limits_{v \in \mathbb{R}} |v g(v)|$. Для $k = 1, \ldots, p$
$$
|d^a_{k}(x, \tau(\theta), \mathbf{U})| = \Biggl | \frac{\partial}{\partial \theta_k} G \Biggl( \frac{x s^{1/2}(\theta, \mathbf{U})}{s^{1/2} (\mathbf{a}, \mathbf{U})} \Biggr) \Biggr | = \frac{|x| U^2_k g(l^a(x, \theta, \mathbf{U}))} {2 s^{1/2} (\theta, \mathbf{U}) s^{1/2} (\mathbf{a}, \mathbf{U})} \leq
$$
\begin{equation}\label{equ_d_k}
\leq g_R (|x| U^2_k) [2 l^a(x, \theta, \mathbf{U}) s^{1/2}(\theta, \mathbf{U}) s^{1/2} (\mathbf{a}, \mathbf{U})]^{-1} = g_R U^2_k [2 s (\theta, \mathbf{U})]^{-1}.
\end{equation}
Последнее выражение можно оценить сверху величиной $2^{-1} g_R U^2_k$. В силу соотношения (\ref{equ_arch_p_inf_sum}) и условия \ref{cond_epsilon_mom_8_plus_e} $\sup \limits_{\mathbf{a} \in \Theta^{\delta}} |y^a_0|^{8 + \beta_0} < \infty$. Следовательно,
первая часть условия ii) выполнена.

Проверка второй части условия ii) технически более сложна. Сначала заметим, что
$$
d^a_{k}(x_1, \tau(\theta^1), \mathbf{U}) - d^a_{k}(x_2, \tau(\theta^2), \mathbf{U}) =
$$
$$
= U^2_k s^{- 1/2} (\mathbf{a}, \mathbf{U}) \Bigl [x_1 g(l^a(x_1, \theta^1, \mathbf{U})) s^{-1/2}(\theta^1, \mathbf{U}) -
$$
$$
- x_2 g( l^a(x_2, \theta^2, \mathbf{U})) s^{-1/2}(\theta^2, \mathbf{U}) \Bigr].
$$
Положим $M(x_0) := \{(x_1, x_2): |x_{1,2}| \leq x_0\}, \; M_a(\alpha) := \{(x_1, x_2, \theta^1, \theta^2): |G(x_1) - G(x_2)| \leq \alpha, \; \| \theta^{1,2} - \mathbf{a} \| \leq  \alpha\}.$ Заметим, что для проверки условия ii) в силу теоремы Лебега о мажорированной сходимости достаточно показать, что при некоторых $C < \infty$ и $\alpha > 0$ для любых $U \in \mathbb{R}^p$ и $n$ выполнено соотношение
\begin{equation}\label{equ_3thm52_1}
\sup \limits_{M_a(\alpha)} |d^a_{k}(x_1, \tau(\theta^1),  \mathbf{U}) - d^a_{k}(x_2, \tau(\theta^2), \mathbf{U})| \leq C U^2_k [s (\mathbf{a}, \mathbf{U})]^{-1},
\end{equation}
а при произвольном фиксированном $\mathbf{U} \in \mathbb{R}^p$ и $\alpha \to 0$
\begin{equation}\label{equ_3thm52_2}
\sup \limits_{M_a(\alpha)} |d^a_{k}(x_1, \tau(\theta^1), \mathbf{U}) - d^a_{k}(x_2, \tau(\theta^2), \mathbf{U})| \to 0
\end{equation}
равномерно по $n$. Соотношение (\ref{equ_3thm52_1}) есть непосредственное следствие (\ref{equ_d_k}). Для доказательства (\ref{equ_3thm52_2}) заметим, что для любого фиксированного $\mathbf{U} \in \mathbb{R}^p$ при $\alpha \to 0$ равномерно по $\mathbf{a} \in \Theta^{\delta}$ и $M_a(\alpha) \cap M(x_0)$ выполняются следующие предельные соотношения:
$$
[x_1 - x_2] \to 0, \quad [s(\theta^1, \mathbf{U}) - s(\theta^2, \mathbf{U})] \to 0,
$$
$$
[l^a(x, \theta^1, \mathbf{U}) - l^a(x, \theta^2, \mathbf{U})] \to 0,
$$
а, следовательно, и $[d^a_{k}(x_1, \tau(\theta^1), \mathbf{U}) - d^a_{k}(x_2, \tau(\theta^2), \mathbf{U})] \to 0$. Кроме того, для достаточно маленького $\alpha > 0$ и фиксированного
$\mathbf{U} \in \mathbb{R}^p$ из условия \ref{cond_distr} следует, что равномерно по $\mathbf{a} \in \Theta^{\delta}$
$$
\sup \limits_{|x| \geq x_0, \| \theta - \mathbf{a} \| \leq \alpha} |x| g(l^a(x, \theta, \mathbf{U})) s^{-1/2}(\theta, \mathbf{U}) \to 0
$$
при $x_0 \to +\infty$. Из последнего соотношения и (\ref{equ_d_k}) следует справедливость второй части условия ii) теоремы \ref{thm_general_proc_op}. Проверка условия iii) проводится
аналогично, это рассуждение мы опускаем. $\Box$

\subsection{Робастность оценки минимального расстояния и GM-оценки.} \label{seq_robustness}

Наряду с асимптотической нормальностью, важным свойством оценки является робастность против отклонений от модели. Мы исследуем поведение MD и GM оценок в схеме засорения
одиночными грубыми выбросами. Напомним ее определение. Пусть $\{z^{\gamma n}_t, \; t = 1, \ldots, n\}$ - независимые бернуллиевские величины с параметром $\gamma$ (т.е. $P(z^{\gamma n}_t = 0) = 1-\gamma, P(z^{\gamma n}_t = 1) = \gamma$), $\{\xi^n_t\}$ - н.о.р. случайные величины с функцией распределения $\mu_{\xi}$ из некоторого фиксированного класса $M_{\xi}$, и последовательности $\{z_t^{\gamma n}\}, \{\xi^n_t\}, \{y^n_t\}$ независимы между собой.  Пусть наблюдаются величины
\begin{equation}\label{equ_outliers}
y_t^{\gamma n} = y^n_t + z^{\gamma n}_t \xi^n_t, \; t = 1 - p, \ldots, 0, 1, \ldots, n.
\end{equation}
В (\ref{equ_outliers}) величины $\{\xi^n_t, \; t = 1, \ldots, n\}$ интерпретируются как грубые выбросы, и описанная схема есть схема засорения данных одиночными грубыми выбросами.

Аналогично случаю выборки, не содержащей засорений, обозначим
$$
\mathbf{Y}^{\gamma n}_{t - 1} := (y^{\gamma n}_{t - 1}, \ldots, y^{\gamma n}_{t - p})^*, \quad \varepsilon^{\gamma n}_t(\theta):=y^{\gamma n}_t / s(\theta, \mathbf{Y}^{\gamma n}_{t-1}), \; t = 1, \ldots, n,
$$
$$
G_n^{\gamma}(x, \theta) := n^{-1} \sum I\{\varepsilon_t^{\gamma n}(\theta) \leq x\},
$$
$$
W^{\gamma}_{n, k}(x, \theta) := n^{-1/2} \sum \varphi_k(\mathbf{Y}^{\gamma n}_{t - 1}, \theta) [I \{ \varepsilon^{\gamma n}_t (\theta) \leq x\} - G^{\gamma}_n(x, \theta)].
$$
Так же, как и ранее, оценки $\mathbf{\hat{a}}_{n, MD}^{\gamma}$ и $\mathbf{\hat{a}}_{n, GM}^{\gamma}$ определяются как решения задач соответственно
$$
K^{\gamma}_{n}(\theta) := \int \|\mathbf{W}^{\gamma}_n(x, \theta)\|^2 d F(x) \to \min \limits_{\theta \in \Theta},
$$
и
$$
\mathbf{l}^{\gamma}_n(\theta) := \sum \varphi(\mathbf{Y}^{\gamma n}_{t-1}, \theta) \psi(\varepsilon_t^{\gamma n}(\theta)) - n^{-1} \sum \varphi(\mathbf{Y}^{\gamma n}_{t-1}, \theta) \sum \psi(\varepsilon_t^{\gamma n}(\theta)) = \mathbf{0}.
$$

Одной из простейших характеристик робастности оценок является их функционал влияния (см., например, \cite{Hampel}, \cite{MY} для определений, а также \cite{Boldin}, \cite{Vyazilov}, где рассчитывается функционал влияния соответственно для ARCH($1$) и GARCH($1, 1$) моделей). Сначала напомним его определение для оценки параметра в произвольной модели. Пусть $\{\mathbf{\hat{a}}_n, n \geq 1\}$ - последовательность оценок параметра $\mathbf{a}$, построенная по выборке с уровнем засорения $\gamma$ (случай $\gamma = 0$ соответствует выборке без засорений). Пусть для любого достаточно малого $\gamma > 0$ при $n \rightarrow \infty$ выполняется соотношение
$$
\mathbf{\hat{a}}_{n} \stackrel{P} \longrightarrow \mathbf{a}^{\gamma}(\mu_{\xi})
$$
для некоторого неслучайного $\mathbf{a}^{\gamma}(\mu_{\xi})$. Тогда функционалом влияния оценки $\mathbf{\hat{a}}_{n}$ называется вектор
$$
\mathbf{IF}(\mathbf{a}^{\gamma}(\mu_{\xi}), \mu_{\xi}) := \lim \limits_{\gamma \rightarrow 0+} \frac{\mathbf{a}^{\gamma}(\mu_{\xi}) - \mathbf{a}}{\gamma}
$$
в том случае, когда она определена. Функционал влияния $\mathbf{IF}(\mathbf{a}^{\gamma}( \mu_{\xi}), \mu_{\xi})$ является главным членом в разложении асимптотического смещения $(\mathbf{a}^{\gamma}(\mu_{\xi}) - \mathbf{a})$ по $\gamma,$
$$
\mathbf{a}^{\gamma}(\mu_{\xi}) - \mathbf{a} = \mathbf{IF}(\mathbf{a}^{\gamma}(\mu_{\xi}), \mu_{\xi}) \gamma + o (\gamma).
$$
Оценка $\mathbf{\hat{a}}_n$ называется робастной в том случае, если $\mathbf{IF}(\mathbf{a}^{\gamma}(\mu_{\xi}), \mu_{\xi})$ определен на $M_{\xi},$ и величина
$$
GES(M_{\xi}, \mathbf{\hat{a}}_n) := \sup \limits_{\mu_{\xi} \in M_{\xi}} \| \mathbf{IF} (\mathbf{a}^{\gamma}(\mu_{\xi}), \mu_{\xi}) \|,
$$
называемая чувствительностью к большим ошибкам, конечна. В этом случае главный член асимптотического смещения $(\mathbf{a}^{\gamma} (\mu_{\xi}) - \mathbf{a})$ равномерно мал при всех возможных засорениях.

Такой подход удобен для исследования оценок, получаемых как решение определенного уравнения, например GM-оценок, так как функционал влияния для них несложно находится (см. \cite{Boldin}). Для MD оценки это сделать сложнее. Поэтому мы предлагаем другую характеризацию. Для общности изложения результат формулируется также и для GM оценки. Следующая теорема является обобщением теоремы \cite[Теорема 2.2]{S} на ARCH($p$) модель, имеющим место для схемы наблюдений с зависящим от $n$ параметром.
\begin{thm}\label{thm_robust}
Пусть выполнены условия \ref{cond_station}, \ref{cond_distr} - \ref{cond_epsilon_mom_8_plus_e}, а также схема наблюдений (\ref{equ_observ_scheme}). Предположим дополнительно, что функция $\varphi$ ограничена
$$
\sup \limits_{\mathbf{U} \in \mathbb{R}^p, \theta \in \Theta} \|\varphi(\mathbf{U}, \theta)\| = L < \infty
$$
и имеют место условия \ref{cond_phi_info_md}, \ref{cond_phi_info_gm}. Тогда для любого $\delta > 0$ существуют $\gamma_0 > 0$ и $N$, такие что для всех $n > N, \gamma < \gamma_0$, любых последовательности $\{\mathbf{a}_n\}$ и распределения выбросов $\mu_{\xi}$
$$
P(\|\mathbf{\hat{a}}_{n, MD}^{\gamma} - \mathbf{a}_n \| > \delta) < \delta, \quad P(\|\mathbf{\hat{a}}_{n, GM}^{\gamma} - \mathbf{a}_n \| > \delta) < \delta
$$
\end{thm}

\textbf{Доказательство.} Мы приведем доказательство лишь для MD оценки. Доказательство в случае GM полностью аналогично. Будем в данном параграфе полагать для краткости $\mathbf{\hat{a}}_n := \mathbf{\hat{a}}_{n, MD}$, $\mathbf{b}^n(x, \theta) := \mathbf{b}^{a_n}(x, \theta)$. Пусть фиксировано некоторое $\delta > 0$. Обозначим $\mathcal{G}$ множество всех вероятностных мер на $\mathbb{R}$. Достаточно показать, что существует $\gamma_0 > 0$, такое, что при $n \rightarrow \infty$
$$
\sup \limits_{\mu_{\xi} \in  \mathcal{G}, \gamma < \gamma_0} P(\| \mathbf{\hat{a}}_{n}^{ \gamma} - \mathbf{a}_n \| > \delta) \longrightarrow 0.
$$
Далее для краткости будем супремум по множеству $\{\mu_{\xi}  \in \mathcal{G}, \gamma < \gamma_0\}$ обозначать просто $\sup.$ Пусть $\delta_0$ - то же, что и в условии \ref{cond_phi_info_md}. Положим
$$
H_n(x, \theta) := \| \mathbf{b}^n(x, \theta) \| - L \gamma_0, \quad d := (2 p + 4) L \gamma_0, \quad c := d + L \gamma_0,
$$
$$
\alpha_0 := c^2 [F(+\infty) - F(-\infty)], \quad \alpha_1 := 4 d^2 [F(+\infty) - F(-\infty)],
$$
$$
\alpha_2 := \delta_0 \delta^2 - \alpha_0 - \delta^{1/2}_0 \delta \alpha_0^{1/2}.
$$
Заметим, что для любого $\theta \in \Theta$, т.ч. $\|\theta - \mathbf{a}_n\| \geq \delta$ и $x \in \mathbb{R}^p$ в силу условия \ref{cond_phi_info_md} и неравенства Коши-Буняковского
$$
\int [H_n(x, \theta) - d]^2 I\{H_n(x, \theta) \geq d\} d F(x) =
$$
$$
= \int [\| \mathbf{b}^n(x, \theta) \| - c]^2 I\{\| \mathbf{b}^n(x, \theta) \| \geq c\}  d F(x) \geq \int \Bigl [ \| \mathbf{b}^n(x, \theta) \|^2 - 2 c \| \mathbf{b}^n(x, \theta) \| \Bigr ] \times
$$
$$
\times I\{\| \mathbf{b}^n(x, \theta) \| \geq c\} d F(x) \geq \int \| \mathbf{b}^n(x, \theta) \|^2 d F(x) - c^2 \int d F(x) -
$$
\begin{equation}\label{equ_rob_alpha_2}
- 2c \int \| \mathbf{b}^n(x, \theta) \| d F(x) \geq  \delta_0 \delta^2 - \alpha_0 - \delta^{1/2}_0 \delta \alpha_0^{1/2} = \alpha_2.
\end{equation}
Несложно проверить, что $\alpha_1 \to 0$, $\alpha \uparrow (\delta_0 \delta^2)$ при $\gamma_0 \to +0$. Следовательно, для достаточно маленького $\gamma_0 > 0$
\begin{equation} \label{equ_rob_1}
\alpha_1 \leq \alpha_2.
\end{equation}
Зафиксируем произвольное такое $\gamma_0$. Покажем, что оно и является искомым. Заметим, что
$$
P(\| \mathbf{\hat{a}}_n^{\gamma} - \mathbf{a}_n \| \leq \delta) \geq P \Bigl (\inf \limits_{\| \theta - \mathbf{a}_n \| \leq \delta} K_n^{\gamma}(\theta) < \inf \limits_{\|\theta - \mathbf{a}_n\| > \delta} K_n^{\gamma}(\theta) \Bigr ) \geq
$$
$$
\geq P \Bigl (K_n^{\gamma}(\mathbf{a}_n) < \inf \limits_{\| \theta - \mathbf{a}_n \| > \delta} K_n^{\gamma}(\theta) \Bigr ).
$$
Таким образом, для завершения доказательства теоремы \ref{thm_robust} достаточно проверить, что
\begin{equation} \label{equ_rob_2}
\sup P \Bigl (K_n^{\gamma}(\mathbf{a}_n) \geq \inf \limits_{\| \theta - \mathbf{a}_n \| > \delta} K_n^{\gamma}(\theta) \Bigr) \to 0.
\end{equation}
В свою очередь, справедливость (\ref{equ_rob_2}) в силу (\ref{equ_rob_1}) обеспечивается следующими двумя соотношениями:
\begin{equation} \label{equ_rob_part_1}
\sup P(K_n^{\gamma}(\mathbf{a}_n) \geq n \alpha_1) \to 0,
\end{equation}
\begin{equation} \label{equ_rob_part_2}
\sup P \Bigl ( \inf \limits_{\| \theta - \mathbf{a}_n \| > \delta} K_n^{\gamma}(\theta) < n \alpha_2 \Bigr ) \to 0.
\end{equation}
Сначала докажем (\ref{equ_rob_part_1}). Очевидно, что $\varepsilon^n_t(\theta) = \varepsilon^{\gamma n}_t(\theta)$ и $\mathbf{Y}^{\gamma n}_{t-1} = \mathbf{Y}^n_{t-1}$, если только $z_t^{\gamma n} = z_{t-1}^{\gamma n} = \ldots = z_{t - p}^{\gamma n} = 0.$ Следовательно, для всех $x \in \mathbb{R}$, $\theta \in \Theta$ и $k = 1, \ldots, p$
$$
|W^{\gamma}_{n, k}(x, \theta) - W_{n, k}(x, \theta)| =
$$
$$
= n^{-1/2} \Biggl | \sum (\varphi_k( \mathbf{Y}^{\gamma n}_{t-1}, \theta) - \varphi_k(\mathbf{Y}^n_{t-1}, \theta)) [I\{\varepsilon^{\gamma n}_t(\theta) \leq x\} - G^{\gamma}_n (x, \theta)] +
$$
$$
+ \sum \varphi_k(\mathbf{Y}^n_{t-1}, \theta) [I\{\varepsilon^{\gamma n}_t(\theta) \leq x\} - I\{\varepsilon^n_t(\theta) \leq x\} - G^{\gamma}_n (x, \theta) + G_n(x, \theta)] \Biggr | \leq
$$
$$ \leq n^{-1/2} \sum |\varphi_k(\mathbf{Y}^{\gamma n}_{t-1}, \theta) - \varphi_k(\mathbf{Y}^n_{t-1}, \theta)|+ \sum \varphi_k(\mathbf{Y}^{\gamma n}_{t-1}, \theta) I\{z_t^{\gamma n} + \ldots + z_{t - p}^{\gamma n} \geq 1\} +
$$
\begin{equation} \label{equ_rob_3}
+ n^{1/2} L |G^{\gamma}_n(x, \theta) - G_n(x, \theta)| \leq (2 p + 4) L n^{-1/2} \sum \limits_{t = 1 - p}^n I\{z_t^{\gamma n} = 1\}.
\end{equation}
Введем следующие последовательности событий:
$$
Q_n := \Bigl \{\forall x \in \mathbb{R}, \| \theta - \mathbf{a}_n \| \geq \delta :\; \| \mathbf{W}_n(x, \theta) \| \geq n^{1/2} H_n(x, \theta) \Bigr \},
$$
$$
R_n := \Bigl \{\sup \limits_{x \in \mathbb{R}, \theta \in \Theta} \| \mathbf{W}^{\gamma}_n(x, \theta) - \mathbf{W}_n(x, \theta) \| < n^{1/2} d \Bigr\}.
$$
В силу соотношений (\ref{equ_wn_exp_wn_diff_estim}), (\ref{equ_rob_3}) и закона больших чисел
\begin{equation} \label{equ_rob_4}
\sup P(\overline{Q}_n \cup \overline{R}_n) \to  0.
\end{equation}
Заметим, что распределение $K_n(\mathbf{a}_n)$ не зависит от $\gamma$ и $\{\xi^n_t\}$ и ${K_n(\mathbf{a}_n) = O_p(1).}$ В силу $(\ref{equ_rob_4})$
$$
\sup P(K_n^{\gamma}(\mathbf{a}_n) \geq n \alpha_1) = \sup P \Biggl (\int \| \mathbf{W}^{\gamma}_n(x, a) \|^2 d F(x) \geq n \alpha_1 \Biggr) \leq
$$
$$
\leq \sup P \Biggl( 2 \Biggl [K_n(\mathbf{a}_n) + \int \| \mathbf{W}^{\gamma}_n(x, \mathbf{a}_n) - \mathbf{W}_n(x, \mathbf{a}_n) \|^2 d F(x)\Biggr] \geq n \alpha_1 \Biggr ) \leq
$$
$$
\leq \sup \Bigl \{ P \Bigl (K_n(\mathbf{a}_n) + n d^2 (F(+\infty) - F(+\infty)) \geq n \alpha_1 / 2 \Bigr) + P(\overline{R}_n ) \Bigr \} \leq
$$
$$
\leq \sup \{ P(K_n(\mathbf{a}_n) / n + \alpha_1 / 4 \geq \alpha_1 / 2) +  P(\overline{R}_n) \} \to 0.
$$
Тем самым (\ref{equ_rob_part_1}) доказано.

Для доказательства (\ref{equ_rob_part_2}) воспользуемся очевидным неравенством:
$$
\| \mathbf{b}_1 + \mathbf{b}_2 \|^2 \geq (\|\mathbf{b}_1\| - \|\mathbf{b}_2\|)^2 I\{\|\mathbf{b}_1\| \geq \|\mathbf{b}_2\| \},
$$
справедливым для любых $\mathbf{b}_1, \mathbf{b}_2 \in \mathbb{R}^p$. Зафиксируем некоторые $x$, $n$ и $\theta,$ т.ч. $\| \theta - \mathbf{a}_n \| > \delta.$ Оценим $\| \mathbf{W}_n^{\gamma}(x, \theta) \|^2$ снизу следующим образом:
$$
\| \mathbf{W}_n^{\gamma}(x, \theta) \|^2 = \| \mathbf{W}_n(x, \theta) + (\mathbf{W}_n^{\gamma}(x, \theta) - \mathbf{W}_n(x, \theta)) \|^2 \geq
$$
$$
\geq \Bigl(\|\mathbf{W}_n(x, \theta) \| - \|\mathbf{W}_n^{\gamma}(x, \theta) - \mathbf{W}_n(x, \theta) \| \Bigr )^2 \times
$$
$$
\times I \Bigl \{\|\mathbf{W}_n(x, \theta) \| \geq \| \mathbf{W}_n^{\gamma}(x, \theta) - \mathbf{W}_n(x, \theta) \| \Bigr \} \geq
$$
$$
\geq \Bigl(\| \mathbf{W}_n(x, \theta) \| - n^{1/2} d \Bigr )^2 I \Bigl \{\| \mathbf{W}_n(x, \theta) \| \geq n^{1/2} d \geq \| \mathbf{W}_n^{\gamma}(x, \theta) - \mathbf{W}_n(x, \theta) \| \Bigr\} \geq
$$
$$
\geq \Bigl( \| \mathbf{W}_n(x, \theta) \| - n^{1/2} d \Bigr)^2 I \Bigl \{ \| \mathbf{W}_n(x, \theta) \| \geq n^{1/2} H_n(x, \theta) \Bigr\} I \Bigl \{H_n(x, \theta) \geq d \Bigr \} I\{R_n\} \geq
$$
$$
\geq n (H_n(x, \theta) - d)^2 I\{Q_n\} I\{R_n\} I\{H_n(x, \theta) \geq d\}.
$$
Следовательно, для любого $\theta,$ т.ч. $\| \theta - \mathbf{a}_n \| > \delta,$ в силу (\ref{equ_rob_alpha_2}) имеем
$$
K_n^{\gamma}(\theta) = \int \| \mathbf{W}_n^{\gamma}(x, \theta) \|^2 d F(x) \geq n I\{Q_n\} I\{R_n\} \times
$$
$$
\times \int (H_n(x, \theta) - d)^2 I\{H_n(x, \theta) \geq d \} d F(x) \geq n I\{Q_n\} I\{R_n\} \alpha_2.
$$
С учетом (\ref{equ_rob_4}) получаем наконец, что
$$
\sup P \Bigl (\inf \limits_{\| \theta - \mathbf{a}_n \| > \delta} K_n^{\gamma}(\theta) < n \alpha_2 \Bigr) \leq
$$
$$
\leq \sup P(n I\{Q_n\} I\{R_n\} \alpha_2 < n \alpha_2) = \sup P(\overline{Q}_n \cup \overline{R}_n) \to 0.
$$
Тем самым (\ref{equ_rob_part_2}), а значит и (\ref{equ_rob_2}), доказаны. $\Box$

\section{Основные результаты.} \label{seq_technical_results}

\subsection{Равномерная оценка коэффициента сильного перемешивания.} \label{seq_regression_uniform_mixing}

Пусть последовательность $\{\varepsilon_t\}$ - та же, что и в (\ref{equ_arch_p}), фиксировано некоторое параметрическое множество $\mathcal{T}$ и для каждого $\tau \in \mathcal{T}$ задан некоторый стационарный процесс $\{\xi^{\tau}_t, \; t \in \mathbb{Z}\}.$ Обозначим для всех $t \in \mathbb{Z}$ $\Xi^{\tau}_t := (\xi^{\tau}_t, \xi^{\tau}_{t-1}, \ldots).$ Мы будем
рассматривать процессы $\{\xi^{\tau}_t\}$, которые для любого $k \in \mathbb{N} \cup \infty$ и $\tau \in \mathcal{T}$  удовлетворяют уравнению
\begin{equation}\label{equ_recurs_equation}
\xi_t^{\tau} = H_k(\varepsilon_t, \ldots, \varepsilon_{t-k+1}, \Xi^{\tau}_{t - k}, {\tau}),
\end{equation}
где $H_1, H_2, \ldots, H_{\infty}-$ некоторые измеримые функции. Несложно проверить, что стационарное решение (\ref{equ_arch_p}) и любого аналогичного уравнения нелинейной авторегрессии (GARCH($p, q$), ARCH($\infty$) и т.д.) в таком виде представляется (для этого достаточно провести несколько итераций уравнения, определяющего соответствующую модель).

Наш первый результат базируется на том, что можно оценить нужный коэффициент сильного перемешивания (в дальнейшем с.п.), зная, насколько быстро убывает зависимость $H_k(\varepsilon_k, \ldots, \varepsilon_1, \mathbf{z}, {\tau})$  от аргумента $\mathbf{z}$ при $k \to \infty$. Точнее, под величиной зависимости $H_k$ от $\mathbf{z}$ здесь мы понимаем расстояние по вариации между с.в. $H_k(\varepsilon_k, \ldots, \varepsilon_1, \mathbf{z}_1, {\tau})$ и $H_k(\varepsilon_k, \ldots, \varepsilon_1, \mathbf{z}_2, {\tau})$ для различных $\mathbf{z}_1$ и $\mathbf{z}_2$. Обозначим $h_{k, t}(\mathbf{z}, {\tau}) := H_k(\varepsilon_t, \ldots, \varepsilon_{t-k+1}, \mathbf{z}, {\tau})$, $d_{TV}$ - расстояние по вариации. Докажем вспомогательную лемму.

\begin{lemma}\label{thm_mixing_coef_estim}
Пусть заданы некоторые польские пространства $\mathcal{X}_1, \; \mathcal{X}_2, \; \mathcal{X}_3$ и независимые случайные элементы $\eta_1$ и $\eta_2$ принимают значения в $\mathcal{X}_1$ и $\mathcal{X}_2$ соответственно, $f(x_1, x_2)$- измеримая функция из $\mathcal{X}_1 \times \mathcal{X}_2$ в $\mathcal{X}_3$. Пусть фиксирован произвольный $x^0_2 \in \mathcal{X}_2$ и $\delta > 0$. Обозначим $P_{\eta_2}$ распределение $\eta_2$, $D_{\delta} := \{x \in \mathcal{X}_2: d_{TV}(f(\eta_1, x), f(\eta_1, x^0_2)) \geq \delta\}$. Тогда коэффициент $\alpha$ с.п. $\sigma$-алгебр $\sigma\{f(\eta_1, \eta_2)\}$ и $\sigma\{\eta_2\}$ удовлетворяет неравенству
$$
\alpha \leq 2 \delta + 2 P(\eta_2 \in D_{\delta}).
$$
\end{lemma}
\textbf{Доказательство.} По определению, $\alpha := \sup |P(f(\eta_1, \eta_2) \in A; \; \eta_2 \in B) - P(f(\eta_1, \eta_2) \in A) P(\eta_2 \in B)|,$ где $\sup$ берется по всем $A \in \mathcal{B} (\mathcal{X}_3), B \in \mathcal{B}(\mathcal{X}_2).$ В силу теоремы Фубини
$$
P(f(\eta_1, \eta_2) \in A; \; \eta_2 \in B) \leq P(f(\eta_1, \eta_2) \in A; \; \eta_2 \in B \backslash D_{\delta}) + P(\eta_2 \in D_{\delta}) =
$$
$$
= \int P(f(\eta_1, x_2) \in A) I\{ x_2 \in B \backslash D_{\delta}\} P_{\eta_2}(d x_2)  + P(\eta_2 \in D_{\delta}) \leq
$$
$$
\leq \int [P(f(\eta_1, x^0_2) \in A) + \delta] I\{ x_2 \in B \backslash D_{\delta}\} P_{\eta_2}(d x_2) + P(\eta_2 \in D_{\delta}) \leq
$$
$$
\leq [ P(f(\eta_1, x^0_2) \in A) + \delta] P( x_2 \in B \backslash D_{\delta}) + P(\eta_2 \in D_{\delta}) \leq
$$
$$
\leq [P(f(\eta_1, \eta_2) \in A) + 2 \delta + P(\eta_2 \in D_{\delta})] P(\eta_2 \in B ) + P(\eta_2 \in D_{\delta}).
$$
Аналогичные оценки дают
$$
P(f(\eta_1, \eta_2) \in A; \; \eta_2 \in B) \geq
$$
$$
\geq [P(f(\eta_1, \eta_2) \in A) - 2 \delta - P(\eta_2 \in D_{\delta})] P(\eta_2 \in B ) - P(\eta_2 \in D_{\delta}).\Box
$$
\begin{thm}\label{thm_expan_recur_mixing_coef}
Пусть процесс $\{\xi^{\tau}_t\}$ удовлетворяет условию (\ref{equ_recurs_equation}), и для некоторого $k_0 \in \mathbb{N}$, борелевской функции $R(\varepsilon_t, \ldots, \varepsilon_{t - k_0 +1}, \mathbf{z}, {\tau})$, $\mathbf{D} \in  \mathcal{B}( \mathbb{R}^{\infty})$, $\mathbf{z}_0 \in \mathbf{D}$ и $\beta \in \mathbb{R}$ и всех $\mathbf{z} \in \mathbf{D}$ имеет место неравенство
$$
d_{TV}(R(\varepsilon_t, \ldots, \varepsilon_{t-k_0+1}, \mathbf{z}_0, {\tau}), R(\varepsilon_t, \ldots, \varepsilon_{t-k_0+1}, \mathbf{z}, {\tau})) \leq \beta.
$$
Обозначим $\alpha$ коэффициент с.п. $\sigma$-алгебр  $\sigma\{\Xi^{\tau}_0\}$ и $\sigma\{R(\varepsilon_{k_0}, \ldots, \varepsilon_1, \Xi^{\tau}_0, {\tau})\}$. Тогда справедливо неравенство $\alpha \leq 2 \beta + 2 P(\Xi^{\tau}_0 \in \overline{\mathbf{D}}).$
\end{thm}
\textbf{Доказательство.} Заметим, что в силу равенства (\ref{equ_recurs_equation}) для $k = \infty$ с.в. $\xi^{\tau}_t$ независима с каждым $\varepsilon_k$ при всех $k > t$. Дальнейшее есть простое применение результата леммы \ref{thm_mixing_coef_estim}. $\Box$

Мы специально выделим как следствие из Теоремы \ref{thm_expan_recur_mixing_coef} случай
$$
R(\varepsilon_t, \ldots, \varepsilon_{t-k+1}, \Xi_{t - k}, {\tau}) = (\xi^{\tau}_{k + p - 1}, \ldots, \xi^{\tau}_k)^*.
$$
Оно понадобится нам при доказательстве результата для ARCH($p$)-модели.
\begin{conseq}\label{thm_expan_recur_finite_coords_mixing_coef}
Предположим, что при $k \in \mathbb{N}$ выполнены условия Теоремы \ref{thm_expan_recur_mixing_coef} для
$$
R_k(\varepsilon_t, \ldots, \varepsilon_{t-k+1}, \mathbf{z}, {\tau}) := (h_{k + p - 1, t + p - 1}(\mathbf{z}, {\tau}), \ldots, h_{k, t}(\mathbf{z}, {\tau}))^*.
$$
Тогда для $\alpha^{\tau}_p(k),$ коэффициента с.п. $\sigma\{\xi^{\tau}_{k + p - 1}, \ldots, \xi^{\tau}_k\}$ и $\sigma\{\Xi^{\tau}_0\},$ при всех $k \in \mathbb{N}$ имеет место неравенство
\begin{equation}\label{equ_finite_coords_mix_coef}
\alpha^{\tau}_p(k) \leq 2 \beta_k + 2 P(\Xi^{\tau}_t \in \overline{\mathbf{D}}_k).
\end{equation}
\end{conseq}
\textbf{Доказательство.} Следует из уравнения (\ref{equ_recurs_equation}), определения $h_{k, t}(\mathbf{z}, \tau)$ и Теоремы \ref{thm_expan_recur_mixing_coef}.

Перейдем к доказательству результата о равномерной оценке коэффициента с.п. для $\{y_t^a\}$. Напомним, что стационарное решение (\ref{equ_arch_p}) удовлетворяет соотношению (разложение Вольтерра, \cite{Giraitis}):
\begin{equation}\label{equ_arch_p_inf_sum}
(y_t^a)^2 = \sum \limits_{l = 0}^{\infty} \sum \limits_{j_1, \ldots, j_l = 1}^{p} a_{j_1} a_{j_2} \ldots a_{j_l} \varepsilon^2_t \varepsilon^2_{t - j_1} \ldots \varepsilon^2_{t -
j_1 - \ldots - j_l}.
\end{equation}
Несложно проверить, что для $k = 1, 2, \ldots$ имеет место следующее представление:
\begin{equation}\label{equ_arch_p_recur_eq}
\begin{array}{ccl}
y_t^a = \varepsilon_t \sigma_{k, t}(\mathbf{Y}_{t - k -1}^a, \mathbf{a}), \\ \sigma^2_{k, t}(\mathbf{z}, \mathbf{a}) := M^0_{t, k}(\mathbf{a}) + M^1_{t, k}(\mathbf{a}) z^2_1 + \ldots + M^p_{t, k}(\mathbf{a}) z_p^2,
\end{array}
\end{equation}
где каждая из случайных величин $M^0_{t, k}(\mathbf{a}), \ldots, M^p_{t, k}(\mathbf{a})$ есть некоторая функция вектора $(\varepsilon_{t-1}, \ldots, \varepsilon_{t-k+1}, \mathbf{a})^*$. Следовательно, имеет место аналогичное (\ref{equ_recurs_equation}) представление $y_t^a = \varepsilon_t \sigma_{k, t}(\mathbf{Y}_{t - k - 1}^a, \mathbf{a}).$ Будем также полагать $\mathbf{Y}^{a}_{k, t}(\mathbf{z}) := (\varepsilon_t \sigma_{k, t}(\mathbf{z}, \mathbf{a}), \ldots, \varepsilon_{t - p + 1} \sigma_{k - p + 1, t - p + 1}(\mathbf{z}, \mathbf{a}))^*.$

Для доказательства мы воспользуемся следствием \ref{thm_expan_recur_finite_coords_mixing_coef} с $\tau = \mathbf{a},$ $\mathbf{z}_0 = (0, 0, \ldots)^*$ и $\mathcal{T} = \Theta^{\delta}$. Основную сложность представляет оценка расстояния по вариации в условии Теоремы \ref{thm_mixing_coef_estim}. Идея ее будет заключаться в том, чтобы представить $\mathbf{Y}^a_{t}$ в виде композиции
$$
\mathbf{Y}^a_{t} = \mathbf{Y}^{a}_{t - p, t} (\mathbf{Y}^{a}_{t - p, t - p}(\mathbf{Y}^{a}_0, \mathbf{a}), \mathbf{a})
$$
и заметить, что при фиксированном $\mathbf{z}$ распределение $\mathbf{Y}^{a}_{t - p, t} (\mathbf{z}, \mathbf{a})$ абсолютно непрерывно, и его плотность непрерывно зависит от $\mathbf{z}$. При этом будет существенно использован мультипликативный вид уравнения (\ref{equ_arch_p}). В свою очередь, $\mathbf{Y}^{a}_{t - p, t - p}(\mathbf{z}, \mathbf{a})$ в некотором смысле "экспоненциально слабо" \, зависит от $\mathbf{z}$ при $k \to \infty$.

Сначала докажем несколько вспомогательных технических лемм.
\begin{lemma}\label{thm_density_of_product}
Пусть $\xi$ и $\eta$ - независимые случайные величины, причем распределение $\xi$ абсолютно непрерывно относительно меры Лебега с плотностью $f_{\xi}(x),$ ограниченной на $\mathbb{R}$. Пусть также $P(|\eta| > d) = 1$ для некоторого $d > 0$. Обозначим $F_{\eta}$ ф.р. $\eta$. Тогда распределение величины $\tau := \xi \eta$ абсолютно непрерывно с плотностью
\begin{equation}\label{equ_prod_density}
f_{\tau}(x) = \int z^{-1} f_{\xi}(x / z) d F_{\eta}(z).
\end{equation}
\end{lemma}
\textbf{Доказательство.} В силу теоремы Фубини $P(\tau \leq x) = \int F_{\xi}(x / z) d F_{\eta}(z)$ для $x \in \mathbb{R}$. По условию, $M := \sup \limits_{x \in \mathbb{R}} f_{\xi}(x) <
\infty.$ Следовательно,
$$
P(\tau \in (x_1, x_2]) = \int [F_{\xi}(x_2 / z) - F_{\xi}(x_1 / z)] d F_{\eta}(z) \leq M (x_2 - x_1) / d.
$$
Отсюда сразу следует абсолютная непрерывность распределения $\tau.$ Сходимость $P(\tau \in (x, x + \delta]) / \delta \to \int z^{-1} f_{\xi}(x / z) d F_{\eta}(z)$ следует из ограниченности
$\delta^{-1} [F_{\xi}(x / z) - F_{\xi}([x + \delta] / z)]$ и теоремы о мажорированной сходимости. $\Box$
\begin{lemma}\label{thm_density_of_product_vector}
Пусть случайный вектор $\mathbf{T}_q(\Xi)$ задается соотношением
$$
\mathbf{T}_q(\Xi) := (\xi_1 \eta_1(\Xi), \xi_2 \eta_2(\xi_1, \Xi), \ldots, \xi_q \eta_q(\xi_1, \ldots, \xi_{q - 1}, \Xi))^*,
$$
где $\xi_1, \ldots, \xi_q$ - н.о.р. случайные величины с распределением, имеющим плотность $f_{\xi}(x),$ $\eta_1, \ldots, \eta_q$ - непрерывные функции, $\Xi$ - случайный элемент, независимый с $\xi_1, \ldots, \xi_q$. Предположим, что соответственно каждой пары с.в. $\eta_j$ и
$\eta_j(\xi_1, \ldots, \xi_{j - 1}, \Xi)$ выполнены условия леммы \ref{thm_density_of_product}. Тогда распределение $\mathbf{T}_q(\Xi)$ абсолютно непрерывно с плотностью
$$
f_T(x_1, \ldots, x_q) := E_{\Xi} \prod \limits_{j = 1}^q f_{\xi}( x_j \eta_j^{-1}(x_1, \ldots, x_j, \Xi) ) \eta_j^{-1}(x_1, \ldots, x_j, \Xi).
$$
\end{lemma}
\textbf{Доказательство.} Доказательство аналогично доказательству леммы \ref{thm_density_of_product}. Положим для краткости $\eta_j(\Xi) := \eta_j(\xi_1, \ldots, \xi_{j - 1}, \Xi).$ По теореме Фубини
$$
P(\xi_j \eta_j(\xi_1, \ldots, \xi_{j - 1}, (\Xi) \in (a_j, b_j], \; j = 1, \ldots, q) =
$$
$$
= P(\xi_j \in (a_j \eta^{-1}_j (\Xi), b_j \eta^{-1}_j(\Xi)], \; j = 1, \ldots, q) =
$$
$$
= E_{\Xi} \int \limits_{a_1 \eta^{-1}_1 (\Xi)}^{b_1 \eta^{-1}_1( \Xi)} f_{\xi}(x_1) \ldots \int  \limits_{a_q \eta^{-1}_q (x_1, \ldots, x_{q - 1}, (\Xi)}^{b_q \eta^{-1}_q(x_1, \ldots, x_{q - 1}, (\Xi)} f_{\xi}(x_q) d x_q \ldots d x_1.
$$
Из такого представления и ограниченности $f_{\xi}$, как и раньше, следует абсолютная
непрерывность распределения $\mathbf{T}_q(\Xi)$. Так как $f_T(x_1, \ldots, x_q)$ непрерывна и равна производной $P(\xi_j \eta_j(\xi_1, \ldots, \xi_{j - 1}, (\Xi) \leq x_j, \; j = 1, \ldots, q)$ по $x_q, \ldots, x_1,$ то она же является искомой плотностью распределения $\mathbf{T}_q(\Xi). \Box$

\begin{conseq}\label{thm_density_of_product_vector_arch}
Распределение вектора $\mathbf{Y}_0^a$ абсолютно непрерывно.
\end{conseq}

Две следующие леммы дают оценку расстояния по вариации между двумя с.в. $\xi \eta_1$ и $\xi \eta_2$, в том случае если с.в. $\eta_1$ и $\eta_2$ близки по вероятности.
\begin{lemma}\label{thm_cmodule_of_density}
Пусть $\xi$, $\eta_1$, $\eta_2$ - с.в., и с. вектора $\xi$ и $(\eta_1, \eta_2)$ независимы. Предположим, что для некоторого $d > 0$ $P(\eta_j > d) = 1, j = 1, 2.$ Пусть распределение $\xi$ абсолютно непрерывно и ее плотность $f_{\xi}(x)$ удовлетворяет условию  \ref{cond_distr}. Положим $K_{\xi} := \int | f'_{\xi}(x) x| d x$, $A := \{|\eta_1 -
\eta_2| \leq \delta\}$. В силу леммы \ref{thm_density_of_product}, существуют плотности $f_{\tau_1}$ и $f_{\tau_2}$ у с.в. $\tau_1 := \xi \eta_1$ и $\tau_2 := \xi \eta_2$. Тогда для них имеет место соотношение
$$
\int |f_{\tau_1}(x) - f_{\tau_2}(x)| d x \leq \delta [K_{\xi} d^{-2} + d^{-1}] + 2 P(\overline{A}).
$$
\end{lemma}
\textbf{Доказательство.} Положим $\eta^A_j := [\eta_j I\{A\} + d I\{\overline{A}\}], \; \tau^A_j = \xi \eta^A_j, \; j = 1, 2.$ Тогда в силу леммы \ref{thm_density_of_product} с.в.
$\tau_j$ и $\tau^A_j$ абсолютно непрерывны. Очевидно, что
$$
|f_{\tau_1}(x) - f_{\tau_2}(x)| \leq |f_{\tau^A_1}(x) - f_{\tau^A_2}(x)| + |f_{\tau^A_1}(x) - f_{\tau_1}(x)| + |f_{\tau^A_2}(x) - f_{\tau_2}(x)|.
$$
и $\int |f_{\tau^A_j}(x) - f_{\tau_j}(x)| dx \leq P(\overline{A}), \; j = 1, 2.$ Поэтому достаточно проверить, что $\int |f_{\tau^A_1}(x) - f_{\tau^A_2}(x)| d x \leq \delta [K_{\xi}
d^{-2} + d^{-1}].$ Заметим, что выражение (\ref{equ_prod_density}) для $f_{\tau^A_j}(x)$ можно переписать следующим образом:$f_{\tau^A_j}(x) = E r(x, \eta^A_j)$, где $r(x, v) := f_{\xi}(x / v) v^{-1}.$ В силу дифференцируемости $f_{\xi}$ (условие \ref{cond_distr}) справедливо неравенство
$$
|f_{\tau^A_1}(x) - f_{\tau^A_2}(x)| \leq E |r(x, \eta^A_1) - r(x, \eta^A_2)| = E |r_v'(x, \eta^A_{1, 2}) [\eta^A_2 - \eta^A_1]|,
$$
для некоторой с.в. $\eta^A_{1, 2}$, т. ч. $\eta^A_{1, 2} \in [\eta^A_1 \wedge \eta^A_2, \eta^A_1 \vee \eta^A_2].$ Следовательно, в силу теоремы Фубини,
$$
\int |f_{\tau^A_1}(x) - f_{\tau^A_2}(x)| d x \leq \int E |r_v'(x, \eta^A_{1, 2}) [\eta^A_2 - \eta^A_1]| dx \leq
$$
$$
\leq \delta E \Biggl \{ \int \Bigl (| f_{\xi}'(x / \eta^A_{1, 2}) x / (\eta^A_{1, 2})^{-3}| + |f_{\xi}(x / \eta^A_{1, 2}) (\eta^A_{1, 2})^{-2}| \Bigr ) d x \Biggl \} =
$$
$$
= \delta E [K_{\xi} (\eta^A_{1, 2})^{-2} + (\eta^A_{1, 2})^{-1}]. \Box
$$
\begin{lemma}\label{thm_cmodule_of_density_vector_conseq}
Пусть выполнены условия леммы \ref{thm_density_of_product_vector} для $(\xi_1, \ldots, \xi_q)^*$ и $\Xi_j, j = 1, 2,$ причем вектор $(\xi_1, \ldots, \xi_q)^*$ независим с $(\Xi_1, \Xi_2)^*$. Положим
$$
A_j := \{|\eta_j(\xi_1, \ldots, \xi_{j - 1}, \Xi_1) - \eta_j(\xi_1, \ldots, \xi_{j - 1}, \Xi_2)| \leq \delta\}, \; j = 1, \ldots, q.
$$
Тогда
$$
d_{TV} (\mathbf{T}(\Xi_1), \mathbf{T}(\Xi_2)) \leq q \delta [K_{\xi} d^{-2} + d^{-1}] + 2 \sum \limits_{j = 1}^q P(\overline{A}_j).
$$
\end{lemma}
\textbf{Доказательство.} Обозначим для краткости $\eta_j(\Xi_i) = \eta_j(\xi_1, \ldots, \xi_{j - 1}, \Xi_i), \; i = 1, 2.$ Аналогично лемме \ref{thm_cmodule_of_density} положим
$$
\eta^A_j(\xi_1, \ldots, \xi_{j - 1}, \Xi_1) := \eta_j(\Xi_1) I\{A_j\} + d I\{\overline{A}_j\}, \; j = 1, \ldots, q.
$$
Пусть случайные векторы $\mathbf{T}^A_1(\Xi_1)$ и $\mathbf{T}^A_2(\Xi_2)$ построены по $\{\eta^A_j(\Xi_i)\}$ вместо $\{\eta_j(\Xi_i)\}$. В силу леммы \ref{thm_density_of_product_vector}, у них существуют плотности соответственно $f_{T^A_i}$ и $f_{T^A_i}$. Как и раньше, $\int |f_{T^A_i}(x_1, \ldots, x_q) - f_{T_i}(x_1, \ldots, x_q)| d \overline{x} \leq \sum \limits_{j = 1}^q P(\overline{A}_j), \; i = 1, 2,$ и поэтому достаточно проверить требуемое неравенство для случая $P(A_j) = 1.$ В обозначениях леммы \ref{thm_cmodule_of_density} имеем:
$$
\int |f_{T^A_1}(x_1, \ldots, x_q) - f_{T^A_2}(x_1, \ldots, x_q)| d \overline{x} =
$$
$$
= \int \Biggl | E \prod \limits_{j = 1}^q r(x_j, \eta^A_j(\Xi_1) ) - E \prod \limits_{j = 1}^q r(x_j, \eta^A_j(\Xi_2)) \Biggr | d \overline{x} \leq
$$
$$
\leq E \int \sum \limits_{j = 1}^q \prod \limits_{i < j} r(x_i, \eta^A_i(\Xi_1)) \prod \limits_{i > j} r(x_i, \eta^A_i(\Xi_2)) |r(x_j, \eta^A_j(\Xi_1) ) - r(x_j, \eta^A_j(\Xi_2))| d \overline{x} =
$$
$$
= E \sum \limits_{j = 1}^q \int \prod \limits_{i < j} r(x_i, \eta^A_i(\Xi_1)) \Biggl ( \int |r(x_j, \eta^A_j(\Xi_1) ) - r(x_j, \eta^A_j(\Xi_2))| d x_j \Biggr ) d x_1 \ldots d x_{j - 1} \leq
$$
$$
\leq \delta ( K_{\xi} d^{-2} + d^{-1} ) E \sum \limits_{j = 1}^q \int \prod \limits_{i < j} r(x_i, \eta_i(\Xi_1)) d x_1 \ldots d x_{j - 1} = q \delta [K_{\xi} d^{-2} + d^{-1}]. \Box
$$
\begin{conseq}\label{thm_expan_recur_mixing_coef_arch}
Предположим, что выполнены условия \ref{cond_station} и \ref{cond_distr}. Тогда существуют $C < \infty$ и $\gamma > 0$, такие что (для всех $\mathbf{a} \in \Theta^{\delta}$) коэффициент с.п. $\sigma\{\mathbf{Y}^a_{t - 1}\}$ и $\sigma\{y^a_{t}, \; t \leq 0\}$ не превосходит $C \exp\{- \gamma t\}.$
\end{conseq}
\textbf{Доказательство.} Фиксируем произвольное $k \in \mathbb{N}.$ Доказательство будет состоять в применении Следствия \ref{thm_expan_recur_finite_coords_mixing_coef}. Для оценки первого члена в неравенстве  (\ref{equ_finite_coords_mix_coef}) применим Лемму \ref{thm_cmodule_of_density_vector_conseq}. Положим, как и раньше, $\mathbf{Y}^{a}_{k, t}(\mathbf{z}) := \varepsilon_t \sigma_{k, t}(\mathbf{z}, \mathbf{a}).$ Несложно проверить, что при $q = p, \; d = 1,$ $\xi_j = \varepsilon_{t - p + j},$
$$
\eta_1(\mathbf{V}) := \sqrt{1 + a_1 V^2_1 + \ldots + a_p V^2_p}, \quad \eta_j(x_1, \ldots, x_{j - 1}, \mathbf{V}) :=
$$
$$
:= \sqrt{1 + a_1 x_{j - 1} \eta_{j - 1}(a) + \ldots + a_{j - 1} x_1 \eta_1(a) + \ldots + a_j V^2_1 + \ldots + a_p V^2_{p - j + 1}}
$$
для $j = 2, \ldots, p$ имеют место равенства $\eta_j(\xi_{t - p + 1}, \ldots, \xi_{t - p + j}, \mathbf{Y}^{\mathbf{a}}_{k - p, t - p}(\mathbf{z})) = \sigma_{k - p + j, t - p + j}(\mathbf{z}, \mathbf{a})$. Следовательно, $\mathbf{T}_p(\mathbf{Y}^{\mathbf{a}}_{k - p, t - p}(\mathbf{z})) = \mathbf{Y}^{\mathbf{a}}_{k, t}(\mathbf{z})$ и выполнены условия Леммы \ref{thm_cmodule_of_density_vector_conseq}. Пользуясь представлениями  (\ref{equ_arch_p_inf_sum}) и (\ref{equ_arch_p_recur_eq}), можно показать, что при некоторых $\gamma_0 > 0, \; C_0 < \infty$
\begin{equation}\label{equ_mj_mom_estim}
E M^j_{0, k}(\mathbf{a}) < C_0 \exp\{-\gamma_0 k\}.
\end{equation}
для всех $\mathbf{a} \in \Theta^{\delta}, \; k \geq 1, \; j = 1, \ldots, p$. Положим $\delta_k := \exp\{ - (\gamma_0 / 2) k\}.$ Тогда в силу неравенства Чебышева и (\ref{equ_mj_mom_estim})
$$
P(\overline{A}_j) = P(|\sigma_{k, t}(\mathbf{z}, \mathbf{a}) - \sigma_{k, t}(\mathbf{0}, \mathbf{a})| > \delta_k) \leq P(|\sigma_{k, t}(\mathbf{z}, \mathbf{a}) - \sigma_{k, t}(\mathbf{0}, \mathbf{a})| \times
$$
$$
\times |\sigma_{k, t}(\mathbf{z}, \mathbf{a}) + \sigma_{k, t}(\mathbf{0}, \mathbf{a})| > 2 \delta_k) = P(|\sigma^2_{k, t}(\mathbf{z}, \mathbf{a}) - \sigma^2_{k, t}(\mathbf{0}, \mathbf{a})| > 2 \delta_k) \leq
$$
$$
\leq P \Biggl ( \sum \limits_{j = 1}^p a_j \varepsilon^2_{t - j} \Biggl [ \sum \limits_{i = 1}^p z_i^2 M^i_{t - j, k - j}(\mathbf{a}) \Biggl ] > 2 \delta_k \Biggr ) \leq 2^{-1} C_0 \delta_k || \mathbf{a}|| E \varepsilon_1^2 ||\mathbf{z}||^2 \leq 2^{-1} C_0 \delta_k ||\mathbf{z}||^2_{2}.
$$
Положим $D_k := \{||\mathbf{z}||^2 \leq 2 C^{-1}_0 \exp\{(\gamma_0 / 4) k\}\}$. Тогда, согласно утверждению Леммы \ref{thm_cmodule_of_density_vector_conseq}, при $\mathbf{z} \in D_k$ расстояние по вариации между $\mathbf{Y}^{ \mathbf{a} }_{k, t}(\mathbf{0})$ и $\mathbf{Y}^{\mathbf{a}}_{k, t} ( \mathbf{ z})$ не превосходит $\beta \exp\{ - (\gamma_0 / 4) k\},$ где $\beta := [K_{\xi} + 2 p + 1].$ Теперь оценим второй член в неравенстве (\ref{equ_finite_coords_mix_coef}). Снова ввиду неравенства Чебышева,
$$
P(\mathbf{Y}^a_0 \in \overline{D}_k) \leq p C_0 E (y_0^a)^2 \exp\{ - \gamma_0 k / 4\} \leq p C_0 E \varepsilon_1^2 (1 - \|\mathbf{b}\|_1 - \delta)^{-1} \exp\{ - \gamma_0 k / 4\}.
$$
Наконец, в силу следствия \ref{thm_expan_recur_finite_coords_mixing_coef}, для $C = 2 \beta + 2 p C_0 E \varepsilon_1^2 (1 - \|\mathbf{b}\|_1 - \delta)^{-1}$ коэффициент с.п. $\sigma\{\mathbf{Y}_t^a\}$ и $\sigma\{y^a_t, \; t \leq 0\}$ не превосходит $ C \exp\{- (\gamma_0 / 4) k\}.\Box$

Рассмотрим два следствия из доказанного результата. Первое из них есть аналог ЗБЧ.
\begin{conseq}\label{thm_LLN}
Пусть для некоторой функции $f$ величины $f(\mathbf{Y}^a_0, \mathbf{a})$ равномерно интегрируемы при $\mathbf{a} \in \Theta^{\delta}$. Тогда
$$
\Bigl | n^{-1} \sum [ f(\mathbf{Y}^a_{t - 1}, \mathbf{a}) - E f(\mathbf{Y}^a_0, \mathbf{a}) ] \Bigr | = o^a_p(1).
$$
\end{conseq}
\textbf{Доказательство.}
Для ограниченной $f$ это утверждение прямо следует из оценки дисперсии, см. \cite[Theorem 5.2]{S}. Для общего случая достаточно ограничить $f$ достаточно большой константой $C$ и воспользоваться равномерной интегрируемостью $f(\mathbf{Y}^n_{0}, \mathbf{a}_n)$. $\Box$

Следующее утверждение является равномерным по $\mathbf{a} \in \Theta^{\delta}$ аналогом Corollary 5.1 из \cite{S}, его доказательство аналогично.
\begin{conseq}\label{thm_sum_proc_max_estim}
Пусть для некоторой функции $f$ $\sup \limits_{\mathbf{a} \in \Theta^{\delta}} E f^4(\mathbf{Y}^a_0, \mathbf{a}) < \infty.$ Тогда
$$
\sup \limits_{x \in \mathbb{R}} \Bigl | n^{-1 / 2} \sum [ f(\mathbf{Y}^a_{t - 1}, \mathbf{a}) I \{\varepsilon_t \leq x\} - G(x) E f(\mathbf{Y}^a_0, \mathbf{a}) ] \Bigr | = O^a_p(1).
$$
\end{conseq}

\textbf{Замечание.} Для обобщения результата следствий \ref{thm_expan_recur_mixing_coef_arch} - \ref{thm_sum_proc_max_estim} на GARCH($p, q$) и аналогичные модели волатильности также можно применить Теорему \ref{thm_expan_recur_mixing_coef}. Пусть интересующая нас модель задается уравнением $\xi^{\theta}_t = \varepsilon_t \sigma(\Xi^{\theta, p}_{t-1}, \Sigma^q_{t - 1}, \theta),$ где $\sigma(\mathbf{z}, \theta)$ - некоторая фиксированная функция,
$$
\Xi^{\theta, p}_{t-1} = (\xi^{\theta}_{t - 1}, \ldots, \xi^{\theta}_{t - p})^*, \quad  \Sigma^q_{t - 1} = (\sigma(\Xi^{\theta, p}_{t-1}, \Sigma^q_{t - 1}, \theta), \ldots, \sigma(\Xi^{\theta, p}_{t - q}, \Sigma^q_{t - q}, \theta))^*.
$$
Тогда в качестве функции $R_k$ из Теоремы \ref{thm_expan_recur_mixing_coef} вместо $(h_{k + q - 1, t + q - 1}(\mathbf{z}, {\theta}), \ldots,$ $h_{k, t}(\mathbf{z}, {\theta}))^*$, как в Следствии \ref{thm_expan_recur_finite_coords_mixing_coef}, надо рассмотреть
$$
(h_{k + p - 1, t + p - 1}(\mathbf{z}, {\theta}), \ldots, h_{k, t}(\mathbf{z}, {\theta}), \sigma_{k + q - 1, t + q - 1}(\mathbf{z}, {\theta}), \ldots, \sigma_{k, t}(\mathbf{z}, {\theta}))^*.
$$
Дальнейшее уточнение выходит за рамки данной статьи и, вероятно, составит тему будущих исследований.

\subsection{Оценка максимума эмпирического процесса.}\label{seq_max_inequality}
В настоящем разделе мы будем рассматривать семейство стационарных случайных процессов $\{\xi^{\tau}_t, \; t \in \mathbb{Z}\}$, где $\tau \in \mathcal{T}$ для некоторого параметрического множества $\mathcal{T}$. Будет предполагаться, что последовательность $\{\varepsilon_t\}$ - та же, что и в (\ref{equ_arch_p}), при каждом $\tau \in \mathcal{T}$ процесс согласован с фильтрацией, образованной $\mathcal{F}_t := \sigma\{\ldots, \varepsilon_0, \varepsilon_1, \ldots, \varepsilon_{t}\}$, и $\xi^{\tau}_t$ не зависит от $\sigma\{\varepsilon_{t+1}, \varepsilon_{t+2}, \ldots\}$. Например, для стационарного решения уравнений типа ARCH и GARCH, или, более общо, для любого процесса вида (\ref{equ_recurs_equation}) это условие выполнено. Будем на протяжении этого параграфа обозначать $\Xi^{\tau}_t := (\xi^{\tau}_t, \xi^{\tau}_{t-1}, \ldots)$ и $\sup := \sup \limits_{\tau \in \mathcal{T}}$. Сначала докажем несколько технических результатов.
\label{seq_auxil_results}
\begin{lemma} \label{thm_m4_max_estim}
Пусть $\eta_1, \eta_2, \ldots, \eta_n$ - с.в., $u_1, u_2, \ldots, u_n$ - набор чисел, причем $E \eta_i^4 < \infty,$ $u_i \geq 0$ при $i = 1, 2, \ldots, n.$ Пусть для некоторого $\alpha > 1$ и произвольных $1 \leq i < j \leq n$ выполнено условие
$$
E (\eta_{i + 1} + \ldots + \eta_j)^4 \leq (u_{i + 1} + \ldots + u_j)^{\alpha}.
$$
Обозначим
$$
S_i := \eta_1 + \ldots + \eta_i; \quad M_n := \max \limits_{0 \leq i \leq n} \min(|S_i|, |S_n - S_i|).
$$
Тогда для некоторого $K,$ не зависящего от $n$ и последовательностей $\{\eta_i\}, \{u_i\}$,
\begin{equation} \label{equ_m4_estim}
E M_n^4 \leq K (u_1 + \ldots + u_n)^{\alpha}.
\end{equation}
\end{lemma}
\textbf{Доказательство.} По индукции. При $n = 1$ соотношение (\ref{equ_m4_estim}) очевидно следует из условия леммы для любого $K \geq 1$. Установим его справедливость для произвольного $n$. Существует такое $j,$ что $u_1 + \ldots + u_{j-1} \leq 1/2,$ $u_{j+1} + \ldots + u_{n} \leq 1/2.$ Положим
$$
M_{1n} := \max \limits_{0 \leq i \leq j - 1} \min(|S_i|, |S_{j - 1} - S_i|); \quad M_{2n} := \max \limits_{j + 1 \leq i \leq n} \min(|S_i - S_j|, |S_n - S_i|).
$$
Тогда, как несложно проверить,
$$
M_n \leq \max(M_{1n} + |S_n - S_{j - 1}|, M_{2n} + |S_j|),
$$
\begin{equation}\label{equ_m4_max_estim}
E M_n^4 \leq E (M_{1n} + |S_n - S_{j - 1}|)^4 + E (M_{2n} + |S_j|)^4.
\end{equation}
Оценим слагаемые в (\ref{equ_m4_max_estim}) по отдельности. Обозначим $u := u_1 + \ldots + u_n.$ Тогда, в силу предположения индукции и условия леммы, $E M_{1n}^4 \leq K(u / 2)^{\alpha}, E (S_n - S_{j - 1})^4 \leq u^{\alpha}.$ Применим неравенство Минковского:
$$
E (M_{1n} + |S_n - S_{j - 1}|)^4 \leq (\|M_{1n}\|_{L^4} + \|S_n - S_{j - 1}\|_{L^4})^4 \leq
$$
$$
\leq ((K (u / 2)^{\alpha})^{1/4} + u^{\alpha / 4})^4 = u^{\alpha} (K^{1/4} 2^{-\alpha/4} + 1)^4.
$$
Аналогично,
$$
E (M_{2n} + |S_j|)^4 \leq u^{\alpha} (K^{1/4} 2^{-\alpha/4} + 1)^4.
$$
Положим $K_{\alpha} = (2^{-1/4} - 2^{-\alpha/4})^{-4}$. Тогда, как несложно проверить, $2 (K_{\alpha}^{1/4} 2^{-\alpha/4} + 1)^4 = K_{\alpha}$. Согласно (\ref{equ_m4_max_estim}),
$$
E M_n^4 \leq 2 u^{\alpha} (K_{\alpha}^{1/4} 2^{-\alpha/4} + 1)^4 = K_{\alpha} u^{\alpha}.
$$
Тем самым, справедливость (\ref{equ_m4_estim}) установлена при $K = K_{\alpha}$ для всех $n$. $\Box$
\begin{conseq}\label{thm_s4_max_estim}
В условиях теоремы \ref{thm_m4_max_estim}
$$
E \max \limits_{i \leq n} S^4_i \leq 2^{\alpha - 1} K_{\alpha} (u_1 + \ldots + u_n)^{\alpha}.
$$
\end{conseq}
\textbf{Доказательство.} Заметим, что для любого $1 \leq i \leq n$ $|S_i| \leq \min(|S_i|, |S_n - S_i|) + |S_n| \leq M_n + |S_n|.$ Следовательно, в силу теоремы \ref{thm_m4_max_estim} и определения$K_{\alpha}$
$$
E \max \limits_{i \leq n} S^4_i \leq E (M_n + |S_n|)^4 \leq ( \|M_n\|_{L^4} + \|S_n\|_{L^4})^4 \leq (K_{\alpha}^{1/4} + 1)^4 (u_1 + \ldots + u_n)^{\alpha} \leq
$$
$$
\leq (K_{\alpha}^{1/4} + 2^{\alpha / 4})^4 (u_1 + \ldots + u_n)^{\alpha} = 2^{\alpha - 1} K_{\alpha} (u_1 + \ldots + u_n)^{\alpha}.\Box
$$
\begin{thm}[Обобщение теоремы 5.1 из \cite{S}]\label{thm_infinit_sum_max_estim}
Пусть фиксированы борелевская функция $f : \mathbb{R}^{\infty} \to \mathbb{R}$, $\mathbf{B} := \{B_j\}_{j=1}^{\infty}$ - разбиение $\mathbb{R}^{\infty},$ т.ч. $B_j \in \mathcal{B}(\mathbb{R}^{\infty})$, и произвольное $\tau \in \mathcal{T}$. Положим $f_j(\mathbf{U}) := f(\mathbf{U}) I\{\mathbf{U} \in B_j\}$ для $\mathbf{U} \in \mathbb{R}^{\infty}$. Пусть выполнены условия:

i) $R(\tau) := \sum \limits_{j = 1}^{\infty} (E f^4_j(\Xi^{\tau}_0))^{1/4} < \infty$.

ii) Для некоторых $N, \; L, \; 1 \leq d < 2$ и произвольных $j \in \mathbb{N}, v \in \mathbb{R}$
$$
\sup \limits_{B_j} |f(\mathbf{U})| \leq N \inf \limits_{B_j} |f(\mathbf{U})|, \quad \sum \limits_{j = 1}^{\infty} \sup \limits_{B_j} |f(\mathbf{U})| I\{\sup \limits_{B_j} |f(\mathbf{U})| \leq v\} \leq L v^d.
$$
Положим $\mathbf{z} := (z_1, z_2, \ldots)^* \in \mathbb{R}^{\infty}$,
\begin{equation} \label{equ_un_def}
S^{\tau}_{n, j}(z) :=  n^{-1/2} \sum  f_j(\Xi^{\tau}_{t - 1}) [I\{\varepsilon_t \leq z\} - G(z)], \quad S^{\tau}_n(\mathbf{z}) := \sum \limits_{j = 1}^{\infty} S^{\tau}_{n, j}(z_j).
\end{equation}
Тогда ряд (\ref{equ_un_def}) для $S^{\tau}_n(\mathbf{z})$ сходится в $L^1$ и п.н. сразу при всех $\mathbf{z} \in \mathbb{R}^{\infty}$, причем
$$
E \sup \limits_{\mathbf{z} \in \mathbb{R}^{\infty}} |S^{\tau}_n(\mathbf{z})| \leq (8 \sqrt{2} K_{3/2})^{1/4} R(\tau) + n^{(d - 2) / 4} 3 M^d L R^d(\tau).
$$
\end{thm}
\textbf{Доказательство.} Сначала установим сходимость ряда (\ref{equ_un_def}) п.н. при каждом $\mathbf{z} \in \mathbb{R}^{\infty}$. Как и при доказательстве теоремы 5.1 из \cite{S}, при помощи неравенства Розенталя (см. \cite[стр. 23]{HH}) устанавливается, что для $z \in \mathbb{R}$
\begin{equation}\label{equ_rosenthal}
E (S^{\tau}_{n, j}(z))^4 \leq C_0 E f^4_j(\Xi^{\tau}_0),
\end{equation}
причем можно положить $C_0 = 8$. Согласно неравенству Минковского и (\ref{equ_rosenthal}),
$$
E |S^{\tau}_{n, j}(z) | \leq ( E (S^{\tau}_{n, j}(z))^4 )^{1/4} \leq C_0^{1/4} (E f^4_j( \Xi^{\tau}_0))^{1/4},
$$
Следовательно, при любых $\mathbf{z} \in \mathbb{R}^{\infty}$ и $n$ ряд (\ref{equ_un_def}) для $S^{\tau}_n(\mathbf{z})$ сходится в $L^1.$ Сходимость п.н. при всех $\mathbf{z}$ очевидна, т.к. только конечное число слагаемых в сумме для $S_n^{\tau}(\mathbf{z})$ отличны от $0$.

Положим $N = N_n := n^2$ и рассмотрим процессы $S^{\tau}_{n, j}(z)$ отдельно на множестве $D_n := \{G^{-1}(i / N_n), i = 0, 1, \ldots, N_n\}$ и на $\mathbb{R} \setminus D_n$. Для завершения доказательства теоремы достаточно проверить, что
\begin{equation}\label{equ_estim_un_BN}
E \sup \limits_{\mathbf{z} \in D^{\infty}_n} |S^{\tau}_n(\mathbf{z})| \leq (8 \sqrt{2} K_{3/2})^{1/4} R(\tau),
\end{equation}
\begin{equation}\label{equ_estim_un_Rinf}
E \sup \limits_{\mathbf{z} \in \mathbb{R}^{\infty}} |S^{\tau}_n(\mathbf{z})| - E \sup \limits_{z \in D^{\infty}_n} |S^{\tau}_n(\mathbf{z})| \leq n^{(d - 2) / 4} 3 N^d L R^d(\tau).
\end{equation}
Аналогично (\ref{equ_rosenthal}) можно показать, что для $1 \leq l \leq m \leq N_n$
$$
E (S^{\tau}_{n, j}(m / N_n) - S^{\tau}_{n, j}(l / N_n))^4 \leq C_0 E f^4_j(\Xi^{\tau}_0) (m / N_n - l / N_n)^{3/2}.
$$
Следовательно, в силу утверждения следствия \ref{thm_s4_max_estim} с $u_i \equiv n N_n^{-1}(E f^4_j(\Xi^{\tau}_0) )^{-3/2}$, $\alpha = 3/2$
$$
E \sup \limits_{z \in D_n} |S^{\tau}_{n, j}(z)| \leq \Bigl (E \sup \limits_{z \in D_n} S^4_{n, j}(z) \Bigr )^{1/4} \leq (2^{1/2} C_0 K_{3/2} E f^4_j(\Xi^{\tau}_0))^{1/4}.
$$
Тем самым установлена справедливость соотношения (\ref{equ_estim_un_BN}). Перейдем к доказательству  (\ref{equ_estim_un_Rinf}). Обозначим $a_0 := -\infty, \; a_i := G^{-1}(i / N_n)$, $A_i := (a_{i-1}, a_i]$ для $i = 1, 2, \ldots, N_n.$ Заметим, что
$$
\sup \limits_{z \in \mathbb{R}} |S^{\tau}_{n, j}(z)| - \sup \limits_{z \in D_n} |S^{\tau}_{n, j}(z)| \leq \max \limits_{i \leq N_n} \sup \limits_{z \in A_i} |S^{\tau}_{n, j}(z) - S^{\tau}_{n, j}(a_{i - 1})|,
$$
$$
\sup \limits_{z \in A_i} |S^{\tau}_{n, j}(z) - S^{\tau}_{n, j}(a_{i-1})| =
$$
$$
= n^{-1/2} \sup \limits_{z \in A_i} \Biggl | \sum f_j(\Xi^{\tau}_{t-1}) [I\{\varepsilon_t \in (a_{i-1}, z]\} - G(z) + G(a_{i-1})] \Biggr | \leq
$$
\begin{equation}\label{equ_estim_un_An-R_diff}
\leq n^{-1/2} \max \limits_{t \leq n} |f_j(\Xi^{\tau}_{t - 1})| \Biggl (\sum I\{\varepsilon_t \in A_i\} + n^{-1}\Biggr).
\end{equation}
Положим $F_n := \max \limits_{t \leq n} |f(\Xi^{\tau}_{t - 1})|$,
$$
j(\mathbf{U}) := \sum \limits_{j = 1}^{\infty} j I\{\mathbf{U} \in B_j\}, \; j_n := \max \limits_{t \leq n} \{j(\Xi^{\tau}_{t - 1})\}, \quad \eta_n := \max \limits_{i \leq N_n} \sum I\{\varepsilon_t \in A_i\} + n^{-1}.
$$
Поскольку $P \Biggl(\sum I\{\varepsilon_t = A_i\} \geq k \Biggr) \leq C_n^k N_n^{-k} \leq N_n^{- k / 2},$ то при $n \geq 2$
$$
P(\eta_n = k ) \leq N_n N_n^{-k / 2} = n^{2 - k}, \quad E \eta^2_n \leq 4 + \sum \limits_{k \geq 3} k^{2} n^{2 - k} \leq 4 + 4 \sum \limits_{k \geq 3} k^{2} 2^{- k} \leq 9.
$$
В силу условия ii) теоремы, (\ref{equ_estim_un_An-R_diff}) и неравенства Гельдера, левая часть  (\ref{equ_estim_un_Rinf}) не превосходит
$$
n^{-1/2} E \eta_n \sum \limits_{j = 1}^{j_n} \sup |f_j| \leq n^{-1/2} E \eta_n \sum \limits_{j = 1}^{\infty} \sup |f_j| I\{\sup |f_j| \leq N F_n\} \leq
$$
$$
\leq n^{-1/2} N^d L E \eta_n |F_n|^d \leq n^{-1/2} N^d L (E \eta^{4/(4 - d)}_n )^{1 - d/4} (E |F_n|^{4})^{d/4} \leq
$$
$$
\leq n^{-1/2 + d / 4} 3 N^d L(E |f(\Xi^{\tau}_0) |^4)^{d/4} \leq n^{(d - 2) / 4} 3 N^d L R^d(\tau).
$$
Тем самым соотношение (\ref{equ_estim_un_Rinf}), а с ним и теорема \ref{thm_infinit_sum_max_estim}, доказаны. $\Box$

\begin{conseq} \label{thm_sup_strict}
Положим
$$
S_n^{0 \tau}(\mathbf{z}) := n^{-1/2} \sum \sum \limits_{j = 1}^{\infty} [f_j(\Xi^{\tau}_{t-1}) I\{\varepsilon_t \leq z_j\} - E f_j(\Xi^{\tau}_0) G(z_j)].
$$
Предположим, что выполнены условия теоремы \ref{thm_infinit_sum_max_estim}. Обозначим $\alpha_{t}({\tau})$ коэффициент с.п. $\sigma$-алгебр $\sigma \{f_j(\Xi^{\tau}_t), \; j \in \mathbb{N}\}$ и $\sigma\{f_j(\Xi^{\tau}_0), \; j \in \mathbb{N}\}$. Предположим, что $Q(\tau ) := \sum \limits_{t = 0}^{\infty} \alpha^{1/2}_{t} ({\tau}) < \infty $. Тогда
$$
E \sup \limits_{\mathbf{z} \in \mathbb{R}^{\infty}} |S_n^{0 \tau}(\mathbf{z})| = R(\tau) [(8 \sqrt{2} K_{3/2})^{1/4} + n^{(d - 2) / 4} 3 M^d L R^{d - 1}(\tau) + (8 Q(\tau))^{1/2}].
$$
\end{conseq}
\textbf{Доказательство.} Положим
$$
s^{\tau}_{j, n}(z) := n^{-1/2} G(z) \sum [f_j(\Xi^{\tau}_{t-1}) - E f_j(\Xi^{\tau}_0)].
$$
Тогда, очевидно, $S_n^{0 \tau}(\mathbf{z}) = S_n^{\tau}(\mathbf{z}) + \sum \limits_{j = 1}^{\infty} s^{\tau}_{j, n}(z_j)$ и имеет место неравенство
\begin{equation}\label{equ_Sno_exp_estim}
E \sup \limits_{\mathbf{z} \in \mathbb{R}^{\infty}} |S_n^0(\mathbf{z})| \leq E \sup \limits_{\mathbf{z} \in \mathbb{R}^{\infty}} |S_n^{\tau}(\mathbf{z})| + E \sum \limits_{j = 1}^{\infty} \sup \limits_{z_j} |s^{\tau}_{j, n}(z_j)|.
\end{equation}
Повторяя для каждого $j$ рассуждения доказательства следствия 5.1 из \cite{S}, получаем, что
$$
E \sup \limits_{z_j \in \mathbb{R}} (s^{\tau}_{j, n}(z_j) )^2  = n^{-1} D \Biggl [ \sum  f_j(\Xi^{\tau}_{t-1}) \Biggr ] \leq
$$
$$
\leq 8 (E f^4_j(\Xi^{\tau}_0))^{1 / 2} \sum \limits_{t = 0}^n \alpha^{1 / 2}_k({\tau}) \leq 8 (E f^4_j(\Xi^{\tau}_0))^{1 / 2} \sum \limits_{k = 0}^{\infty} \alpha^{1 / 2}_k({\tau}).
$$
Следовательно,
$$
E \sum \limits_{j = 1}^{\infty} \sup \limits_{z_j} |s^{\tau}_{j, n}(z_j)| \leq  \sum \limits_{j = 1}^{\infty} [E \sup \limits_{z_j} |s^{\tau}_{j, n}(z_j)|^2]^{1/2} \leq
$$
$$
\leq R(\tau) \Biggl [ 8 \sum \limits_{k = 0}^{\infty} \alpha^{1 / 2}_k({\tau}) \Biggr ]^{1/2} = R(\tau) [8 Q(\tau)]^{1/2}.
$$
В силу утверждения теоремы \ref{thm_infinit_sum_max_estim} и (\ref{equ_Sno_exp_estim}), следствие доказано$.\Box$

\begin{lemma}\label{thm_2vars_mom_estim}
Пусть $\eta$ - случайная величина на $(\Omega, \mathcal{F}, P)$, такая что для некоторого $\beta > 0$ $M := E \eta^{4 + \beta} < \infty.$ Пусть $\{A_j, \; j = 1, 2, \ldots\}, A_j \in \mathcal{F}$ образуют разбиение $\Omega.$ Предположим, что $L := \sum \limits_{j=1}^{\infty} (P(A_j))^{\beta/4 (4 + \beta)} < \infty.$ Тогда
$$
\sum \limits_{j=1}^{\infty} (E \eta^4 I_{A_j})^{1/4} \leq M^{4/(4 + \beta)} L.
$$
\end{lemma}
\textbf{Доказательство.} Согласно неравенству Гельдера,
$$
E \eta^4 I_{A_j} \leq (E \eta^{4 + \beta})^{4/(4 + \beta)} (P(A_j))^{\beta/(4 + \beta)}.
$$
Следовательно,
$$
\sum \limits_{j=1}^{\infty} (E \eta^4 I_{A_j})^{1/4} \leq (E \eta^{4 + \beta})^{1/(4 + \beta)} \sum \limits_{j=1}^{\infty} (P(A_j))^{\beta/4(4 + \beta)} = M^{1/(4 + \beta)} L.\Box
$$

Следующая теорема является основным результатом раздела. Она устанавливает оценку максимума эмпирического процесса $T_n$, построенного по наблюдениям из стационарного процесса $\xi_t^{\tau}$ вида (\ref{equ_recur_proc}). Пусть фиксированы компакт $\mathbf{Z} \subset \mathbb{R}^q$ и, для каждого $\tau \in \mathbf{Z}$, борелевские функции $\Delta^{\tau}(x, \mathbf{z}, \mathbf{U}) : \mathbb{R} \times \mathbf{Z} \times \mathbb{R}^{\infty} \rightarrow \mathbb{R}$, $\lambda^{\tau}(\mathbf{z}, \mathbf{U}) : \mathbf{Z} \times  \mathbb{R}^{\infty} \rightarrow \mathbb{R}$.
\begin{thm}\label{thm_general_proc_op}
Пусть семейство процессов $\{\xi^{\tau}_t\}$ удовлетворяет условиям следствия \ref{thm_sup_strict}, где в качестве $\alpha_t(\tau)$ выступает коэффициент с.п. $\sigma\{( \lambda(\tau, \Xi^{\tau}_t), \; \Delta^{\tau}(x, \mathbf{z}, \Xi^{\tau}_t)), \; \mathbf{z} \in \mathbf{Z}, x \in \mathbb{R}\}$ и $\sigma\{\Xi^{\tau}_0\}$. Положим для $\mathbf{z} \in \mathbf{Z}, \; x \in \mathbb{R}$
$$
T^{\tau}_n(x, \mathbf{z}) := n^{-1/2} \sum \lambda^{\tau}(\mathbf{z}, \Xi^{\tau}_{t-1}) I\{\varepsilon_t \leq \Delta^{\tau}(x, \mathbf{z}, \Xi^{\tau}_{t-1})\} - n^{1/2} E \lambda^{\tau}(\mathbf{z}, \Xi^{\tau}_0) G(\Delta^{\tau}(x, \mathbf{z}, \Xi^{\tau}_0)).
$$
Пусть для некоторых $\beta > 0$ и $\mathbf{z}^0 \in \mathbf{Z}$ выполнены следующие условия:

i) Для всех $\tau \in \mathbf{Z}$, $\mathbf{U} \in \mathbb{R}^{\infty}$ и $x \in \mathbb{R}$ $\Delta^{\tau}(x, \mathbf{z}^0, \mathbf{U}) = x.$

ii) Для всех $\tau \in \mathbf{Z}$, $\mathbf{U} \in \mathbb{R}^{\infty}$ и $k = 1, \ldots, q$ функция $G(\Delta^{\tau}(x, \mathbf{z}, \mathbf{U}))$ имеет частные производные $d^{\tau}_k(x, \mathbf{z}, \mathbf{U})$ по $z_k$, причем
$$
\sup E (d^{\tau}_k(x, \mathbf{z}^0, \Xi^{\tau}_0))^{4 + \beta} < \infty, \quad \sup E | d^{\tau}_k(x^1, \mathbf{z}^1, \Xi^{\tau}_0) - d^{\tau}_k(x^2, \mathbf{z}^2,  \Xi^{\tau}_0) |^{4 + \beta} \to 0.
$$
при $\| \mathbf{z}^1 - \mathbf{z}^2 \| \to 0, \; |G(x^1) - G(x^2)| \to 0.$

iii) Для $\tau \in \mathbf{Z}$, $k = 1, \ldots, q, \; \mathbf{z} \in \mathbf{Z}$ и всех $\mathbf{U} \in \mathbb{R}^{\infty}$ функция $\lambda^{\tau}(\mathbf{z}, \mathbf{U})$ имеет частные производные $h^{\tau}_k(\mathbf{z}, \mathbf{U})$ по $z_k$, причем для некоторой $H$, т.ч. $E |H(\Xi^{\tau}_0)|^2 < \infty,$
$$
\sup E |\lambda^{\tau}(\mathbf{z}^0, \Xi_0^{\tau})|^2 < \infty, \quad |h^{\tau}_k(\mathbf{z}, \mathbf{U})| \leq H(\mathbf{U}), \quad \sup E | h^{\tau}_k(\mathbf{z}^1, \Xi^{\tau}_0) - h^{\tau}_k(\mathbf{z}^2, \Xi^{\tau}_0) |^2 \to 0,
$$
при $\| \mathbf{z}^1 - \mathbf{z}^2 \| \to 0.$

Тогда для любого $\rho > 0$
\begin{equation}\label{equ_estim_vn}
\sup \limits_{\mathbf{z} \in \mathbf{Z}, x \in \mathbb{R}} (T^{\tau}_n(x, \mathbf{z}) \pm n^{1/2} \rho ||\mathbf{z} - \mathbf{z}^0||)^{\mp} = O^{\tau}_p(1).
\end{equation}
\end{thm}
\textbf{Доказательство.} Сначала покажем, что для произвольного фиксированного $\rho > 0$ существует $\delta > 0$, т.ч. (\ref{equ_estim_vn}) справедливо с заменой $\sup \limits_{\mathbf{z} \in \mathbf{Z}, x \in \mathbb{R}} \{\cdot\}$ на $\sup \limits_{ \mathbf{Z}_\delta, \; x \in \mathbb{R}} \{\cdot\},$ где $\mathbf{Z}_{\delta} = \mathbf{Z} \cap \{\mathbf{z}: \|\mathbf{z} - \mathbf{z}^0\|_1 \leq \delta\}$. Без ограничения общности будем считать, что $\mathbf{z}^0 = \mathbf{0}$ и $\mathbf{Z} \subset \mathbb{R}^q_+.$ Доказательство для общего случая аналогично. Будем опускать верхний индекс $\tau$ в записи $\Delta^{\tau}, \; d^{\tau}_k, \; \lambda^{\tau}_k, \; h^{\tau}_k$.

\textbf{Доказательство для случая $\mathbf{z} \in \mathbf{Z}_{\delta}$.}
Зафиксируем некоторое $N \in \mathbb{N}$. Пусть $x_l := G^{-1}(l / N),$ $l = 0, 1, \ldots, N,$ $A_l := (x_{l - 1}, x_l]$. Для $\mathbf{U} \in \mathbb{R}^{\infty}, \; k = 1, \ldots, q, x \in A_l, \; D > 0$ положим
$$
d_{k, \delta}^+(x, \mathbf{U}) := \sup \limits_{\mathbf{z} \in \mathbf{Z}_{\delta}, v \in A_l} d_{k}(v, \mathbf{z}, \mathbf{U}); \quad d_{k, \delta}^-(x, \mathbf{U}) := \inf \limits_{\mathbf{z} \in \mathbf{Z}_{\delta}, v \in A_l} d_{k}(v, \mathbf{z}, \mathbf{U}),
$$
$$
h_{k, \delta}^+(\mathbf{U}) := I\{H(\mathbf{U}) \leq D\} \sup \limits_{\mathbf{z} \in \mathbf{Z}_{\delta}} h_{k}(\mathbf{z}, \mathbf{U}); \quad h_{k, \delta}^-(\mathbf{U}) := I\{H(\mathbf{U}) \leq D\} \inf \limits_{\mathbf{z} \in \mathbf{Z}_{\delta}} h_{k}(\mathbf{z}, \mathbf{U}).
$$
Пусть $\mathbf{B}^{+} := \{B^+_j, j \in \mathbb{N}\}, \mathbf{B}^{-} := \{B^-_j, j \in \mathbb{N}\}$ - некоторые разбиения $\mathbb{R}^{\infty}$, которые будут выбраны позднее. Определим два вспомогательных процесса:
$$
T_{n, \delta}^{\pm}(x, \mathbf{z}) := n^{-1/2} \sum \lambda^{\pm}_{\delta} (\mathbf{z}, \Xi^{\tau}_{t-1}) I\{\varepsilon_t \leq \Delta_{\delta}^{\pm}(x, \mathbf{z}, \Xi^{\tau}_{t-1})\} -
$$
$$
- n^{1/2} E \lambda^{ \mp}_{\delta}(\mathbf{z}, \Xi^{\tau}_0) G(\Delta^{ \mp}_{\delta}(x, \mathbf{z}, \Xi^{\tau}_0)),
$$
где
$$
\lambda^{\pm}_{\delta}(\mathbf{z}, \mathbf{U}) := \lambda(\mathbf{0}, \mathbf{U}) + \sum \limits_{k = 1}^q z_k h_{k , \delta}^{\pm}(\mathbf{U}),
$$
$$
d^+_{k, \delta, j}(x) := \sup \limits_{\mathbf{U} \in B^+_j} d^+_{k, \delta}(x, \mathbf{U}),  \quad d^-_{k, \delta, j}(x) := \inf \limits_{\mathbf{U} \in B^-_j} d^-_{k, \delta}(x, \mathbf{U}),
$$
$$
\Delta^{\pm}_{\delta}(x, \mathbf{z}, \mathbf{U}) := G^{-1} \Biggl (G(x) + \sum \limits_{k = 1}^q z_k \sum \limits_{j = 1}^{\infty} I\{\mathbf{U} \in B^{\pm}_j\} d^{\pm}_{k, \delta, j}(x) \Biggr )
$$
и $G^{-1}$ доопределена как $G^{-1}(x) = - \infty$, $x \leq 0$, $G^{-1}(x) = \infty$, $x \geq 1.$ Для краткости обозначим
$$
\lambda^{\pm}_{t} := \lambda^{\pm}_{\delta}(\mathbf{z}, \Xi^{\tau}_{t-1}), \quad v^{\pm}_t := I\{\varepsilon_t \leq \Delta_{\delta}^{\pm}(x, \mathbf{z}, \Xi^{\tau}_{t-1})\}.
$$
Положим $S^H_n(D) := n^{-1} \sum I\{ H( \Xi^{\tau}_{t-1}) \geq D\}.$ Заметим, что
$$
T_{n, \delta}^-(x, \mathbf{z}) - n^{1/2} \|\mathbf{z}\|_1 S^H_n(D) \leq T_n(x, \mathbf{z}) \leq T_{n, \delta}^+(x, \mathbf{z}) + n^{1/2} \|\mathbf{z}\|_1  S^H_n(D).
$$
Легко видеть, что величина $n^{1/2} \sup E S^H_n(D)$ не зависит от $x$ и $\mathbf{z}$ и, в силу условия iii) и леммы \ref{thm_LLN}, может быть сделана сколь угодно малой при достаточно большом $D$. Следовательно, справедливость (\ref{equ_estim_vn}) с $\mathbf{z} \in \mathbf{Z}_{ \delta}$ будет обеспечена, если мы покажем, что для любого $\rho > 0$ существуют такие $\delta > 0, \; N \in \mathbb{N}$ и разбиения $\mathbf{B}^{\pm}$, что для всех $l \leq N$
\begin{equation}\label{equ_vn_pm_op}
\sup \limits_{\mathbf{z} \in \mathbf{Z}_{ \delta}, x \in A_l} \Bigl |T_{n, \delta}^{\pm}(x, \mathbf{z}) - n^{1/2} E [\lambda^{\pm}_1 v^{\pm}_1 - \lambda^{\mp}_1 v^{\mp}_1] \Bigr | = O^{\tau}_p(1),
\end{equation}
и для всех $\mathbf{z} \in \mathbf{Z}_{ \delta}, \; x \in \mathbb{R}$
\begin{equation}\label{equ_cov-e2_less rho}
E \lambda^+_1 v^+_1 - E \lambda^-_1 v^-_1 \leq \rho ||\mathbf{z}||.
\end{equation}
Мы докажем (\ref{equ_vn_pm_op}) для $T_{n, \delta}^+(x, \mathbf{z})$, для $T_{n, \delta}^-(x, \mathbf{z})$ доказательство полностью аналогично. Будем в дальнейшем предполагать, что $1 \leq l \leq N$ фиксировано и $x \in A_l$.

Разбиение $\mathbf{B}^+$ будем выбирать следующим образом. Зафиксируем некоторое $s > 0,$ и положим $b_j := (\exp\{s (j - 1)\} - 1)$ для $j \geq 1$. Пусть $j \mapsto (m_1(j), \ldots, m_q(j))$ - такая нумерация $\mathbb{N}^q,$ что $\max \limits_{1\leq i \leq m} (m_i(j))$ не убывает. Положим для $j = 1, 2, \ldots$
$$
B^{+}_j := \{\mathbf{U} \in \mathbb{R}^{\infty}: |d^{+}_{k, \delta}(x_l, \mathbf{U})| \in [b_{m_k(j)}, b_{m_{k}(j)+1}), k = 1, \ldots, q\}.
$$
Фиксируем произвольные $\mathbf{U}_j \in B_j$, тогда для  всех $\mathbf{U} \in B^+_j$ $\Delta^+_{\delta}(x, \mathbf{z}, \mathbf{U}) = \Delta^+_{\delta}(x, \mathbf{z}, \mathbf{U}_j)$. Положим
$$
\Delta^+_{\delta, j}(x, \mathbf{z}) := \Delta^+_{\delta}(x, \mathbf{z}, \mathbf{U}_j), \quad \mathbf{r} := \Biggl ( \Delta^+_{\delta, 1}(x, \mathbf{z}), \Delta^+_{\delta, 2}(x, \mathbf{z}), \ldots \Biggr )^* \in \mathbb{R}^{\infty},
$$
$$
f^0(\mathbf{U}) = \lambda(0, \mathbf{U}), \quad z_0 \equiv 1, \quad f^k(\mathbf{U}) = h_{k , \delta}^+(\mathbf{U}), \quad k = 1, \ldots, q.
$$
Теперь в обозначениях теоремы \ref{thm_infinit_sum_max_estim} можно записать:
$$
T_{n, \delta}^+(x, \mathbf{z}) - n^{1/2} [E \lambda^+_1 v^+_1 - E \lambda^-_1 E v^-_1] = n^{-1/2} \sum \limits_{k = 0}^q z_k \sum \limits_{j = 1}^{\infty} \sum I\{\Xi^{\tau}_{t-1} \in B^+_j \} \times
$$
$$
\times[f^k(\Xi^{\tau}_{t-1}) I\{\varepsilon_t \leq r_j\}- E f^k(\Xi^{\tau}_0) G(r_j))] = S_{0, n}^{0 \tau} (\mathbf{r}) + \sum \limits_{k = 1}^q z_k S_{k, n}^{0 \tau} (\mathbf{r}),
$$
Проверим, что при $k = 0, 1, \ldots, q$ для $f^k$ выполняются условия i) и ii) теоремы \ref{thm_infinit_sum_max_estim}. Заметим, что
\begin{equation}\label{equ_fk_sum_estim}
\sum \limits_{j = 0}^{\infty} ( E [f^k_j(\Xi^{\tau}_0)]^4)^{1/4} = \sum \limits_{j = 0}^{\infty} \Bigl( E [f^k(\Xi^{\tau}_0)]^4 I\{\Xi^{\tau}_0 \in B^+_j\} \Bigr)^{1/4}.
\end{equation}
Для оценки правой части (\ref{equ_fk_sum_estim}) воспользуемся леммой \ref{thm_2vars_mom_estim}. Проверим, что выполнено ее условие с $A_j = \{\Xi^{\tau}_0 \in B_j^+\}, \; \eta = d^+_{k, \delta}(x_j, \Xi^{\tau}_0)$. Положим $M_k := E [d^+_{k, \delta} (x_l, \Xi^{\tau}_0)]^{4 + \delta}, M := M_1 + \ldots + M_q,$
$$
m_{max}(j) := \max \limits_{k \leq q} m_k(j), \quad k_{max}(j) := \mathop{argmax} \limits_{k \leq q} m_k(j).
$$
Заметим, что согласно неравенству Чебышева
$$
P(\Xi^{\tau}_0 \in B^+_j) \leq P(|d^{+}_{k_{max}(j), \delta}(x_l, \Xi^{\tau}_0)| \geq b_{m_{max}(j)}) \leq [b_{m_{max}(j)}]^{-4 - \delta} M_{k_{max}(j)},
$$
и, согласно определению $\{b_j\}$ и взаимной однозначности отображения $m: \mathbb{N} \to \mathbb{N}^q$,
$$
\sum \limits_{j = 1}^{\infty} [P(\Xi^{\tau}_0 \in B^+_j)]^{\delta/(4 + \delta)} \leq M^{\delta/(4 + \delta)} \sum \limits_{j=1}^{\infty} [b_{m_{max}(j)}]^{-\delta} \leq
$$
$$
\leq M^{\delta/(4 + \delta)} \sum \limits_{m=1}^{\infty} (m^q - (m - 1)^q) b^{-\delta}_m \leq
$$
$$
\leq M^{\delta /(4 + \delta)} \sum \limits_{m=1}^{\infty} (m^q - (m - 1)^q) \exp(- \delta s m) < \infty.
$$
Итак, в силу леммы \ref{thm_2vars_mom_estim} условие i) теоремы \ref{thm_infinit_sum_max_estim} выполнено. Проверим, что имеет место ii) (с $N = \exp\{s\}$). С учетом конкретного вида множеств $\{B^+_j\}$ достаточно проверить, что для некоторых $d < 2, \; L < \infty$ и любых $l \in \mathbb{N},$ $k = 1, \ldots, q$
$$
\sum \limits_{m = 1}^l (m^q - (m - 1)^q) b_m \leq L b_l^d.
$$
Но, действительно, при всех $l \in \mathbb{N}$
$$
\sum \limits_{m = 1}^l (m^q - (m - 1)^q) b_m \leq \sum \limits_{m = 1}^l (m^q - (m - 1)^q) b_l \leq l^q b_l \leq L b^d_l,
$$
для произвольного $d > 1, \; L := \max \limits_{l \in \mathbb{N}} \{ l^{- q} \exp(s l (d - 1)) \}.$ Итак, условия теоремы \ref{thm_infinit_sum_max_estim} выполнены (требуемое условие на коэффициент с.п. выполнено автоматически, так как каждая из функций $I\{\mathbf{U} \in \mathbf{B}^{+}_j\} f^k(\mathbf{U})$ выражается через $\lambda$ и $\Delta$). В силу Теоремы \ref{thm_infinit_sum_max_estim} имеем
$$
\sup \limits_{\|\mathbf{z}\| \leq \delta, x \in \mathbb{R}} |T_{n,  \delta}^+(x, \mathbf{z}) - n^{1/2} E [\lambda^+_1 v^+_1 - \lambda^-_1 v^-_1]| \leq
$$
$$
\leq \sup \limits_{\mathbf{r} \in \mathbb{R}^{\infty}} |S_{0, n}^{0 \tau}(\mathbf{r})| + \delta \sum \limits_{k = 1}^q \sup \limits_{\mathbf{r} \in \mathbb{R}^{\infty}} | S_{k, n}^{0 \tau} (\mathbf{r})| = O^{\tau}_p(1).
$$
Тем самым справедливость (\ref{equ_vn_pm_op}) доказана. Перейдем к соотношению (\ref{equ_cov-e2_less rho}).

Преобразуем левую часть (\ref{equ_cov-e2_less rho}) следующим образом:
\begin{equation}\label{equ_cov-e2_lambda&v}
|E \lambda^+_1 v^+_1 - E \lambda^-_1 v^-_1| \leq E |\lambda^+_1| [v^+_1 - v^-_1] + E |\lambda^+_1 - \lambda^-_1| v^-_1.
\end{equation}
Второе слагаемое в (\ref{equ_cov-e2_lambda&v}) ограничено сверху величиной $2 E [\lambda^+_1 - \lambda^-_1],$ что, в свою очередь, не превосходит $2 \|\mathbf{z}\|_1 \sum \limits_{k = 1}^q E [h_{k, \delta}^+(\Xi^{\tau}_0) - h_{k, \delta}^-(\Xi^{\tau}_0)]$. Согласно условию
iii), последняя сумма может быть сделана сколь угодно маленькой соответствующем выбором $\delta.$ Оценим теперь первое слагаемое. Заметим, что для $j = 1, 2, \ldots$, $k = 1, \ldots, q$ и любых $x \in A_l,$ $\mathbf{U}_j^{\pm} \in B^{\pm}_j$ выполнено
соотношение $d^{\pm}_{k, \delta}(x, \mathbf{U}_j^{\pm}) \in [b_{m_k(j)}, b_{m_k(j) + 1}]$. Следовательно, при тех же $j, k, x, \mathbf{U}^{\pm}$ верно, что $0 \leq \pm (d^{\pm}_{k, \delta, j}(x_l) - d^{\pm}_{k, \delta}(x, \mathbf{U}^{\pm})) \leq (b_{m_k(j) + 1} - b_{m_k(j)})$. Имеем:
$$
E |\lambda^+_1| [v^+_1 - v^-_1] = E |\lambda^+_1| [G(\Delta_{\delta}^+(x_l, \mathbf{z}, \Xi^{\tau}_0)) - G(\Delta_{\delta}^-(x_l, \mathbf{z}, \Xi^{\tau}_0))] \leq
$$
$$
\leq E |\lambda^+_1| \sum \limits_{j = 1}^{\infty} \Biggl [ I\{\Xi^{\tau}_0 \in B^+_j\} \sum \limits_{k = 1}^q z_k d^+_{k, \delta, j}(x_l) - I\{\Xi^{\tau}_0 \in B^-_j\} \sum \limits_{k = 1}^q z_k d^-_{k, \delta, j}(x_l) \Biggr ] \leq
$$
$$
\leq \| \mathbf{z} \|_{\infty} E |\lambda^+_1| \Biggl ( \sum \limits_{j = 1}^{\infty} ( I\{ \Xi^{\tau}_0 \in B^+_j\} + I\{\Xi^{\tau}_0 \in B^-_j\}) \sum \limits_{k = 1}^q [b_{m_k(j) + 1} - b_{m_k(j)} ] +
$$
$$
+ \sum \limits_{k = 1}^q [d^+_{k, \delta}(x_l, \Xi^{\tau}_0) - d^-_{k, \delta}(x_l, \Xi^{\tau}_0)] \Biggr ) \leq
$$
\begin{equation}\label{equ_pm_expect_diff_estim}
\begin{array}{ccl}
\leq \| \mathbf{z} \|_{\infty} E |\lambda^+_1| \Bigl [ (e^s - 1)( |d^+_{k, \delta}(x_l, \Xi^{\tau}_0)| + \\ + |d^-_{k, \delta}(x_l, \Xi^{\tau}_0)| ) + [d^+_{k, \delta}(x_l, \Xi^{\tau}_0) - d^-_{k, \delta}(x_l, \Xi^{\tau}_0)] \Bigr ].
\end{array}
\end{equation}
Для доказательства (\ref{equ_cov-e2_less rho}) достаточно показать, что коэффициент при $\|\mathbf{z}\|_{\infty}$ в (\ref{equ_pm_expect_diff_estim}) стремится к $0$ при $s, \delta \rightarrow 0.$ Это несложно проверить, используя условие ii).

\textbf{Доказательство для общего случая.} Пусть фиксировано некоторое $\rho > 0$. Тогда, как было показано выше, для некоторого $\delta$ соотношение (\ref{equ_estim_vn}) выполнено c $\mathbf{z} \in \mathbf{Z}_\delta$. Положим $\gamma := \delta \rho / 4.$ Положим $\Lambda(\mathbf{U}) := | \lambda(\mathbf{0}, \mathbf{U})| + \mathop{diam} (\mathbf{Z})  H(\mathbf{U}),$ тогда для всех $\mathbf{z} \in \mathbf{Z}$ и $\mathbf{U} \in \mathbb{R}^{\infty}$ $|\lambda(\mathbf{z}, \mathbf{U})| \leq \Lambda(\mathbf{U})$. Существуют $\alpha > 0$, $N$ и разбиение $\mathbf{Z} = \mathbf{Z}_1 \sqcup \ldots \sqcup \mathbf{Z}_N$, $\mathbf{z}_1 \in \mathbf{Z}_1, \ldots, \mathbf{z}_N \in \mathbf{Z}_N,$ а также $-\infty = x_0 < x_1 < \ldots < x_N = \infty,$ такие что для любых $l, m \leq N$
$$
E \sup \limits_{A_l \times \mathbf{Z}_m} |\lambda(\mathbf{z}, \Xi^{\tau}_0) - \lambda(\mathbf{z}_m, \Xi^{\tau}_0)| < \gamma,
$$
$$
E \Lambda(\Xi^{\tau}_0) \sup \limits_{A_l \times \mathbf{Z}_m} |G(\Delta(x, \mathbf{z},
\Xi^{\tau}_0)) - G(\Delta(x_l, \mathbf{z}_m, \Xi^{\tau}_0))| < \gamma,
$$
где $A_l := (x_{l-1}, x_l].$ Положим при любых $1 \leq m, l \leq N, \; \mathbf{U} \in \mathbb{R}^{\infty}$ и всех $x \in A_l, \; \mathbf{z} \in \mathbf{Z}_m$
$$
\lambda^+(\mathbf{z}, \mathbf{U}) := \sup \limits_{\mathbf{z}_+ \in \mathbf{Z}_m} \lambda(\mathbf{z}_+, \mathbf{U}), \quad \lambda^-(\mathbf{z}, \mathbf{U}) := \inf \limits_{\mathbf{z}_- \in \mathbf{Z}_m} \lambda(\mathbf{z}_-, \mathbf{U}),
$$
$$
\Delta^+(x, \mathbf{z}, \mathbf{U}) := \sup \limits_{A_l \times \mathbf{Z}_m} \Delta(x_+, \mathbf{z}_+, \mathbf{U}), \quad \Delta^-(x, \mathbf{z}, \mathbf{U}) := \inf \limits_{A_l \times \mathbf{Z}_m} \Delta(x_-, \mathbf{z}_-, \mathbf{U}).
$$
$$
T^{\pm}_n(x, \mathbf{z}) := n^{-1/2} \sum \lambda^{\pm}(\mathbf{z}, \Xi^{\tau}_{t-1}) I\{\varepsilon_t \leq \Delta^{\pm}(x, \mathbf{z}, \Xi^{\tau}_{t-1})\} -
$$
$$
- n^{1/2} E \lambda^{\mp}(\mathbf{z}, \Xi^{\tau}_0) G(\Delta^{\mp}(x, \mathbf{z}, \Xi^{\tau}_0)).
$$
Для краткости примем также $\lambda^{\pm}_m(U) = \lambda^{\pm}(\mathbf{z}^{0, m}, U), \Delta^{\pm}_{m, l}(U) = \Delta^{\pm}(x_l, \mathbf{z}^{0, m}, U).$ Несложно показать, что для
любых $1 \leq m, l \leq N$
\begin{equation}\label{equ_vn_pm_op_2}
T^{\pm}_n(x_l, \mathbf{z}^{0, m}) \mp n^{1/2} E [\lambda^+_m (\Xi^{\tau}_0) G(\Delta^+_{m, l}(\Xi^{\tau}_0)) - \lambda^-_m( \Xi^{\tau}_0) G(\Delta^-_{m, l}(\Xi^{\tau}_0))] = O^{\tau}_p(1).
\end{equation}
Заметим, что согласно определению $\{x_l\}$ и $\{\mathbf{z}_m\}$
$$
E [\lambda^+_m (\Xi^{\tau}_0) G(\Delta^+_{m, l}(\Xi^{\tau}_0)) - \lambda^-_m( \Xi^{\tau}_0) G(\Delta^-_{m, l}(\Xi^{\tau}_0))] =
$$
$$
= E \lambda^+_m(\Xi^{\tau}_0) [G(\Delta^{+}_{m, l}(\Xi^{\tau}_0)) - G(\Delta^{-}_{m, l}(\Xi^{\tau}_0))] +
$$
$$
+ E[\lambda^{+}_m(\Xi^{\tau}_0) - \lambda^{-}_m(\Xi^{\tau}_0)] G(\Delta^{-}_{m, l}(\Xi^{\tau}_0)) \leq 2 \gamma + 2 \gamma = 4 \gamma.
$$
Следовательно, в силу (\ref{equ_vn_pm_op_2}) и определения $\lambda^{\pm}_m$ и $\Delta^{\pm}_{m, l}$,
$$
\sup \limits_{\mathbf{z} \in \mathbf{Z} \backslash \mathbf{Z}_{\delta}, x \in \mathbb{R}} (T_n(x, \mathbf{z}) \pm n^{1/2} \rho ||\mathbf{z}||)^{\mp} \leq \sup \limits_{\mathbf{z} \in \mathbf{Z} \backslash \mathbf{Z}_{\delta}, x \in \mathbb{R}} (T^{\mp}_n(x, \mathbf{z}) \pm n^{1/2} \rho \delta)^{\mp} \leq
$$
$$
\leq \max \limits_{m, l \leq N} \sup \limits_{A_l \times \mathbf{Z}_m} (T^{\mp}_n(x_l, \mathbf{z}^m) \pm n^{1/2} \rho \delta)^{\mp} \leq \max \limits_{m, l \leq N} n^{1/2} [E \lambda^{\pm}_m(\Xi^{\tau}_0) G(\Delta^{\pm}_{m, l}(\Xi^{\tau}_0)) -
$$
$$
- E \lambda^{\mp}_m(\Xi^{\tau}_0) G(\Delta^{\mp}_{m, l}(\Xi^{\tau}_0)) \pm \rho \delta]^{\mp} + O^{\tau}_p(1) = O^{\tau}_p(1). \Box
$$

\end{document}